\documentclass[a4paper, 11pt]{article}
\usepackage{bbm}
\usepackage{amsfonts}
\usepackage{amsmath,amssymb,amsthm,amsfonts}
\usepackage{amssymb}
\usepackage{graphicx}
\usepackage{epstopdf}
\usepackage{latexsym}
\usepackage{amsmath}
\usepackage{amsthm}
\usepackage{pdfsync}
\usepackage[colorlinks=true]{hyperref}

\hypersetup{dvipdfmx,
       colorlinks,
       linkcolor=blue,
       filecolor=red,
       urlcolor=blue,
       citecolor=red,
       pdftitle={Vortex structures},
       pdfauthor={ Jun YANG, Shenzhen University},
       pdfsubject={   PDE         },
       pdfkeywords={   Geometric flows},
       pdfproducer={ps2pdf}}
\makeatletter
\def\@currentlabel{2.1}\label{e:dispaa}
\def\@currentlabel{2.21}\label{e:dispau}
\def\@currentlabel{2.22}\label{e:dispav}
\def\@currentlabel{2.23}\label{e:dispaw}
\def\@currentlabel{2.24}\label{e:dispax}
\def\theequation{\thesection.\@arabic\c@equation}

\makeatother
\def\theequation{\thesection.\@arabic\c@equation}
\makeatother

\renewcommand{\theequation}{\thesection.\arabic{equation}}
\newtheorem{theorem}{Theorem}[section]
\newtheorem{corollary}[theorem]{Corollary}
\newtheorem{proposition}[theorem]{Proposition}
\newtheorem{remark}[theorem]{Remark}
\newtheorem{lemma}[theorem]{Lemma}

\newtheorem{definition}{Definition}
\newcommand{\bremark}{\begin{remark} \em}
\newcommand{\eremark}{\end{remark} }

\newcommand{\R}{{\mathbb R}}
\newcommand{\ve}{\varepsilon}
\newcommand{\F}[2]{\frac{#1}{#2}}

\newcommand{\BE}{\begin{equation}}
\newcommand{\BEN}{\begin{equation*}}
\newcommand{\EE}{\end{equation}}
\newcommand{\EEN}{\end{equation*}}
\newcommand{\BL}{\begin{lemma}}
\newcommand{\EL}{\end{lemma}}
\newcommand{\BT}{\begin{theorem}}
\newcommand{\ET}{\end{theorem}}
\newcommand{\BP}{\begin{proposition}}
\newcommand{\EP}{\end{proposition}}
\newcommand{\BC}{\begin{corollary}}
\newcommand{\EC}{\end{corollary}}

\allowdisplaybreaks

\begin{document}


\title{\bf Geometric solitons with vortex structures for some geometric flows from Minkowski spaces: part I}
\author{
Youde Wang
\\
Academy of Mathematics and Systems Science,
Chinese Academy of Sciences,
\\
Beijing 100190, P. R. China.
wyd@math.ac.cn
\\
\\
Jun Yang\\  College of Mathematics and Computational Sciences,
Shenzhen University,
\\
Nanhai Ave 3688, Shenzhen 518060,  P. R. China.
jyang@szu.edu.cn
}




\date{}
\maketitle
\begin{abstract}
By a perturbation approach, we construct geometric solitons with various vortex structures(vortex pairs, vortex rings)
 for some geometric flows(Wave maps, Shr\"odinger flows) from Minkowski spaces to ${\mathbb S}^2\subset\R^3$.
\end{abstract}

\setcounter{equation}{0}
\section{Introduction}\label{section1}
The well-known Landau-Lifshitz (LL) equation is a fascinating nonlinear evolution equation,
which describes the dynamics of several important physical systems such as ferromagnets,
moving space curves, etc. and has intimate connections with many of the well known integrable soliton equations,
including nonlinear Schr\"{o}dinger  equations and modified KdV equations.
It admits rich dynamical structures including solitons, dromions, vortices, etc.

\medskip
In many situations ferromagnetic materials may be viewed as continuous media.
For instance, all the spins are locked together by the strong exchange interaction
for a single-domain ferromagnetic nanoparticle at sufficiently low temperatures.
In order to solve the so-called momentum problem in ferromagnets \cite{H},
G. Volovik proposed a simple model of delocalized electrons in \cite{Vo}.
The model introduces the normal velocity of the fermionic liquid as an additional hydrodynamical variable,
describing the background fermionic vacuum,
which is inspired by superfluid motion at $T =0$ in $\text{He}^3$--A (\cite{HoM, HoM1}).
The purpose of this approach is the restoration of the correct linear momentum density
by using these "hydrodynamical" variables.
Moreover, for planar magnetic systems the existence of topologically nontrivial vortex
configurations requires the fluid to be rotational with non-vanishing vorticity function \cite{MPS}.

\medskip
To be precise, the density and the normal velocity of the fermionic liquid,
introduced in \cite{Vo}, are denoted by $\rho$ and $\text{v}$ respectively.
Now, since the vorticity of $\text{v}$ is proportional to the density $\rho$ of
the topological current of the magnetic configuration,
one may formulate a simple model of an anisotropic ferromagnetic fluid,
which arises as a modification of the phenomenological Landau-Lifshitz equation,
with the aim of trying to solve the momentum problem.
A magnetic fluid, or a spin-liquid, is a continuum of material particles,
each of which is endowed with a magnetic moment.
In the continuum approximation the spatial variations of the physical quantities occur
over a large number of particles, so that the medium is characterized by the macroscopic field $m(r, t)$,
i.e., the local magnetization with the constraint $|m|^ 2 = 1$,
the velocity field $\text{v}(r, t)$ in the Lagrange picture,
and the density of particles $\rho(r, t)$.

\medskip
In a local moving frame with velocity $\text{v}(r, t)=(v_1,\cdots, v_n)$,
the precession motion of $m=(m_1, m_2, m_3)$ is described by the Landau-Lifshitz equation in $n$-dimensions
\begin{equation}\label{MLL}
\partial_tm + v^j\partial_jm = -m\times\nabla E_e(m)
\end{equation}
and $v_i$ is a sort of vector gauge potential,
whose tensor field, the vorticity in the hydrodynamical language, is given by
\begin{equation}\label{MLL1}
\partial_iv_j-\partial_jv_i=2m\cdot(\partial_im \times \partial_jm).
\end{equation}
Here $E_e(S)$ denotes the effective energy, which is of the general form
$$
E_e(m)=\frac{1}{2}\int_{\mathbb{R}^n}(|\nabla m|_g^2 + \sum_{i=1}^3\beta_im_i^2)\,\mathrm{d}x
+ \int_{\mathbb{R}^n}m\cdot H \,\mathrm{d}x,
$$
where the $\beta_i$'s are the anisotropy constants,
$H$ is the external magnetic field, and the inner product in $\mathbb{R}^n$
is defined by the constant diagonal metric $g=\text{diag}(\alpha_1,\cdots,\alpha_n)$
where $\alpha_i$ are the exchange coefficients. By the definition of metric $g$,
we have
$$
|(x_1,\cdots,x_n)|^2_g=\sum_{i=1}^n\alpha_ix_i^2,\quad\quad v^j\partial_jm
=\sum_{i=1}^ng^{ij}v_i\partial_jm, \quad\quad |\nabla m|_g^2 =\sum_{i=1}^n\alpha_i|\partial_im|^2.
$$
In the sequel, we will use the summation convection.
For compressible fluid we need, in addition, to take account of the continuity equation
$$\partial_t\rho=\text{div}(\rho\text{v}), \quad\quad \text{v}=(v_1, \cdots, v_n),$$
and to assume that the exchange and anisotropy coefficients in $E_e$ depend on $\rho$.
However, in the present paper we restrict ourselves to the incompressible case,
that is $\rho$ is a constant (see \cite{MPS}).

\medskip
In the case of a ferromagnetic boundary condition
$$m(r,t)\rightarrow (0,0,-1)\quad\quad \text{for}\quad\quad \sum_{i=1}^nx_i^2\rightarrow\infty,$$
the $n$-dimensional domain space can be compactified as a topological $n$-sphere $\mathbb{S}^n$, and $m$ defines mappings $\mathbb{S}^n$ into $\mathbb{S}^2$ corresponding to homotopy classes belonging to $\pi_n(\mathbb{S}^2)$. In the case of $2$ and $3$ dimensions, the mapping class $[m]$ can be characterized by the following topological index:

(1) for the case $n=2$, it is the integer-valued topological degree of the mapping from $\mathbb{S}^2$ onto $\mathbb{S}^2$ (topological charge)
$$Q=\frac{1}{4\pi}\int_{\mathbb{R}^2}m\cdot(\partial_1m \times \partial_2m)\,\mathrm{d}x;$$

(2) for the case $n=3$, it is just the Hopf invariant of $m$ from $\mathbb{S}^3$ into $\mathbb{S}^2$ (linking number)
$$L_H=-\frac{1}{(8\pi)^2}\int_{\mathbb{R}^3}\epsilon_{ijk}v_im\cdot(\partial_jm \times \partial_km)\,\mathrm{d}x,$$
where the Levi-Civita symbol for $n=3$ is given by
$$\epsilon_{ijk}=\frac{(i-j)(j-k)(k-i)}{2}.$$
Moreover, the ``magnetic vorticity" (topological current) is defined by
$$J_i=\frac{1}{2}\epsilon_{ijk}m\cdot(\partial_jm \times \partial_km),$$
which is a conserved quantity in the case of $3$ dimensions and satisfies $\partial^iJ_i=0$.

\medskip
In the case of one dimension the above model is just the usual Landau-Lifshitz equation\cite{LL}.
The Ishimori model \cite{Ishimori} can be recovered only for particular choices of the
metric $g=\text{diag}(1, -1)$. It is well known that the Ishimori model is the first example of integrable
classical spin model in $2+1$ dimensions.

\medskip
According to the signature of the metric $g$ we have the different
cases $(+,\cdots, +)$ and $(+,\cdots,+, -,\cdots, -)$, possibly after a rescaling
of the spatial variables.
One usually denotes $(\mathbb{R}^n, g)$ by $\mathbb{R}^{K,N}$ with $K+N=n$,
and $\mathbb{R}^{1,N}$ is just the Minkowski space.
In mathematics, we need to consider the following Landau-Lifshitz equation
with hydrodynamical term (see \cite{KN})
\begin{eqnarray}\label{ML1}
\begin{array}{lll}
\partial_tm &=& m\times \left((\partial^2_1+\cdots+\partial^2_K -\partial^2_{K+1}-\cdots-\partial^2_{K+N})m - \sum_{i=1}^3 \beta_i m_i\vec{e}_i-H\right)\\
& & - b\left(\sum_{j=1}^Kv_j\partial_jm - \sum_{j=K+1}^{K+N}v_j\partial_jm\right),
\end{array}
\end{eqnarray}
\begin{equation}\label{ML2}
\partial_iv_j-\partial_jv_i=2m\cdot(\partial_im \times \partial_jm),
\end{equation}
with the incompressible condition
$$\text{div}(\text{v})=0,$$
where $b$ is a real constant, $\vec{e}_1=(1, 0, 0)$, $\vec{e}_2=(0, 1, 0)$ and $\vec{e}_3=(0, 0, 1)$. This equation unifies almost all Landau-Lifshitz models.

\medskip
In the present work we will be concerned only with $b=0$ and the isotropic higher dimension case,
in the absence of an external field.
Our task turn into  considering the following problems from  Minkowski space
into a standard unit sphere ${\mathbb S}^2$ in $\mathbb{R}^3$
\begin{equation}\label{Ishimoriproblem}
\partial_{t}m\,=\,m\times\big(\Box m \,+\,|Dm|^2m\big)
\quad \mbox{in} \quad \R\times\R^{1,N}
\quad\mbox{with}\quad m(t,\tau, s)\in\mathbb{S}^2\subset\R^3,
\end{equation}
and the equation for its stationary solutions
\begin{equation}\label{Wavemaps}
\Box m \,+\,|Dm|^2m\,=\,0
\quad
\mbox{in} \quad \R^{1,N}
\quad\mbox{with}\quad m(\tau, s)\in\mathbb{S}^2\subset\R^3.\\
\end{equation}
Here, the differential operators are defined by
\begin{align*}
\begin{aligned}
\Box m\,\equiv\,
\triangle_\tau m\,-\,\triangle_s m
\,=\,
\frac{\partial^2 m}{\partial \tau^2}
\,-\,
\sum_{j=1}^{N}\frac{\partial^2 m}{\partial s_j^2},\qquad\quad
\\
|Dm|^2\,\equiv\,
\big|\nabla_\tau m\big|^2\,-\,\big|\nabla_s m\big|^2
\,=\,
\Big|\frac{\partial m}{\partial \tau}\Big|^2
\,-\,\sum_{j=1}^{N}\Big|\frac{\partial m}{\partial s_j}\Big|^2.
\end{aligned}
\end{align*}
Note that for the convenience of future use of notation,
we here denote $\tau$ and $s$ the variables in Minkowski space.

\medskip
It is easy to see that  (\ref{Ishimoriproblem}) and (\ref{Wavemaps}) are some types of geometric flows.
More precisely (\ref{Ishimoriproblem}) can be regarded as the Schr\"odinger flows for
maps from a Minkowski space into ${\mathbb S}^2$ in the sense of \cite{D1} (also see \cite{SongWang}).
While (\ref{Wavemaps}) can be viewed as a generalization of the ${\mathbb S}^2$
$\sigma$-model which is a two dimensional scalar field equation with ${\mathbb S}^2$
as the target manifold, i.e., $S:\mathbb{R}^{1,2}\rightarrow {\mathbb S}^2$.
Usually, one calls (\ref{Wavemaps}) for $S: \mathbb{R}^{1,N} \rightarrow {\mathbb S}^2$ as the wave map equation.

\medskip
From the viewpoint of analysis, many mathematicians made contributions to the Landau-lifshitz equation without hydrodynamical term. For the Schr\"odinger flows from the Euclidean space $\mathbb{R}^n$ into $\mathbb{S}^2$,
despite some serious efforts (see e.g. \cite{BSS,BIKT}, \cite{GTZ,GH}, \cite{GKT, GKT2},
\cite{IK1,IK2}, \cite{S2,S3}) and the references therein),
some basic mathematical issues such as global well-posedness and
global in time asymptotic for the equation (\ref{Ishimoriproblem}) remain unknown.
One also studied the Schr\"odinger map flow from $\mathbb{R}^n$ or a compact Riemannian
manifold into a K\"ahler manifold (\cite{S1}, \cite{D1}-\cite{D4}, \cite{RRS} and the references therein).
If one is interested in one-dimensional wave (plane-wave) solutions of (\ref{Ishimoriproblem}),
that is, $m: \R \times \R \to {\mathbb S}^2$ (or ${\mathbb S}^1$),
a lot were known as (\ref{Ishimoriproblem}) becomes basically an integrable system (see \cite{GTZ}).
The problem in 2-D or higher dimensions are much more subtle. Even though it is possible
to obtain weak solutions of (\ref{Ishimoriproblem})
(see \cite{D1}-\cite{D2}, \cite{DG1, DG2}, \cite{GTZ}, \cite{M}, \cite{V}),
one does not know if such weak solutions are classical (smooth) or unique.

\medskip
One also studied the Cauchy problems of the Ishimori model and some existence results were established recently. We refer to \cite{KN} for more details.
Especially, C. Kenig, G. Ponce and L. Vega in \cite{KenigPonceVega} have ever studied the following Schr\"odinger equation which is analogous to the Schr\"odinger flow from Lorentzian manifold:
\begin{align}
\left\{
\begin{aligned}
&\frac{\partial u}{\partial t}\,=\,i\,\mathfrak{L} u\,+\,P(u,\nabla u,{\bar u},\nabla{\bar u}),
\\
&u(0,x)\,=\,u_0(x),
\end{aligned}\right.
\end{align}
where $u=u(t, x)$ is a complex valued function from $\R\times\R^N$,
$\mathfrak{\Phi}$ is a non-degenerate second-order operator
$$
\mathfrak{L}\,=\,\sum_{j\leq K}\partial^2_j\,-\,\sum_{j>K}\partial^2_j,
$$
for some $K\in\{1,\cdots, n\}$, and $P:{\mathbb C}^{2n+2}\rightarrow{\mathbb C}$ is a polynomial satisfying certain constraints. They proved the local well-posedness of the above initial value problem in appropriate Sobolev spaces. No doubt, it is a hard task to settle the existence and uniqueness problems of such a class of Cauchy problems (\ref{Ishimoriproblem}).

\medskip
It is well-known that vortex dynamics is a natural paradigm for the field of chaotic motion
and modern dynamical system theory.
The vortex dynamics provides some physically profound examples of Hamiltonian systems of
infinite dimensions, attracting much interest in connection with chaotic phenomena in dynamical systems.
From the viewpoint of fluid mechanism, it is of important significance to study the behavior of the vortex
solutions to such Ishimori type equation.
Indeed, some exact vortex solutions for the isotropic case ($\beta=0$) are considered in the case $n\geq2$.
For instance, in \cite{GP,MPS} the two dimensional case without external magnetic field was studied
and a wide class of solutions can be generated by using time-dependent gauge transformations.
For the Landau-Lifshitz equation (\ref{MLL})
from a standard Euclidean space $\mathbb{R}^n$ with $b=0$, $\beta=0$ and $H=0$,
by our knowledge there is no existence result on vortex solutions in analytic aspects.
However, one has established some existence results on vortex solutions to
the equation with $b=0$, $\beta=(0, 0, 1)$ and $H=0$ by employing the Lyapunov-Schmidt
reduction method (\cite{linwei,weiyang}).

\medskip
A natural problem is whether or not one can employ an applicable analysis method
to construct some solutions to (\ref{Ishimoriproblem}), which are of vortex structure. This is the main aim of the present paper.
More precisely, we intend to construct some geometric soliton (or solitary wave)
solutions with various vortex structures to the above problem (\ref{Ishimoriproblem})
and (\ref{Wavemaps}).
Thus, let's to recall some well-known facts and some known results on the existence of vortex
solutions with finite energy to the Ishimori model on $\mathbb{R}^{1,1}\times\mathbb{R}$
and the Landau-Lifshitz equation (\ref{ML1}) from the Euclidean space $\mathbb{R}^n$ with $b=0$
and the anisotropic higher dimension case, in the absence of an external field.

\subsection{Anti-holomorphic Ishimori model}
In the case $n=2$, $b=1$ and $\beta=(\beta_1,\beta_2,\beta_3)=0$, the above Landau-Lifshitz equation (\ref{ML1})-(\ref{ML2})  without external magnetic field for $m: \mathbb{R}^{1,1} \times\mathbb{R}\rightarrow S^2$ reads
\begin{equation}\label{DE}
\partial_t m + v_1\partial_1m- v_2\partial_2m = m\times(\partial^2_1-\partial^2_2)m,
\end{equation}
\begin{equation}\label{DE1}
\partial_1v_2-\partial_2v_1=2m\cdot(\partial_1m\times\partial_2m),
\end{equation}
with the incompressible condition
\begin{equation}\label{DE2}
\partial_1v_1 + \partial_2v_2 = 0.
\end{equation}
The following conservation law holds true
$$\partial_0\mathfrak{J}_0-\partial_1\mathfrak{J}_1 + \partial_2\mathfrak{J}_2 = 0,$$
where
$$\mathfrak{J}_0= |\partial_1m|^2+|\partial_2m|^2,$$
$$\mathfrak{J}_1=v_1\mathfrak{J}_0-2\partial_1m\cdot(m\times(\partial^2_1-\partial^2_2)m)+ 2m\cdot(\partial_1m\times\partial_2^2m-\partial_1\partial_2 m\times\partial_2m),$$
and
$$\mathfrak{J}_2=2\partial_2m\cdot(m\times(\partial^2_1-\partial^2_2)m) + v_2\mathfrak{J}_0 -2m\cdot(\partial_1^2m\times\partial_1\partial_2m-\partial_1 m\times\partial_2 m).$$

The energy functional
$$E(m)=\int_{\mathbb{R}^2}(|\partial_1m|^2+|\partial_2m|^2) \,\mathrm{d}x$$
and the topological charge are conserved quantities, which are related by the Bogomolny type inequality
$$E(m)\geq 8\pi|Q|.$$
It follows from the following evident inequality
$$\int_{\mathbb{R}^2}|\partial_im\pm\epsilon_{ij}m\times\partial_jm|^2\,\mathrm{d}x\ge 0,$$
and is saturated by time dependent spin configurations satisfying the self-duality equations
$$\partial_im\pm\epsilon_{ij}m\times\partial_jm = 0.$$
If we consider the standard unit sphere $\mathbb{S}^2$ as the Riemann sphere for a complex plane, we have the following stereographic projections
$$m_1+im_2=\frac{2\zeta}{1+|\zeta|^2}, \quad\quad m_3=\frac{1-|\zeta|^2}{1+|\zeta|^2},$$
where $\zeta(x,y,t)$ is a complex value function. In stereographic projection form, the self-dual equation corresponding the first sign can be rewritten
$$\zeta_{\bar{z}}(x,y,t)=\partial_{\bar{z}}\zeta(x,y,t)=0$$
while for the second sign we have the anti-holomorphicity
$$\zeta_{z}(x,y,t)=\partial_{z}\zeta(x,y,t)=0.$$
Here $z=x+iy$ and $\partial_z=\frac{1}{2}(\partial_x-i\partial_y)$, $\partial_{\bar{z}}=\frac{1}{2}(\partial_x-i\partial_y)$. These conditions written in terms of the real $\text{Re}\,\zeta$ and imaginary $\text{Im}\,\zeta$ parts of function $\zeta$, representing the Cauchy-Riemann equations, describe the incompressible and irrotational fluid flow with the velocity potential $\text{Re}\,\zeta$ and the stream function $\text{Im}\,\zeta$.

The equation (\ref{DE2}) can be solved in terms of the stream function of the flow, $v_1=-\partial_2\varphi$, $v_2=\partial_1\varphi$ so that we get the so called Ishimori Model
$$\partial_tm = m\times(\partial_1^2m-\partial_2^2m) + \partial_2\varphi\partial_1m+\partial_1\varphi\partial_2m,$$
$$\Delta\varphi=(\partial_1^2 + \partial_2^2)\varphi=2m\cdot(\partial_1m\times\partial_2m).$$
In terms of complex variables,
$$v_+ = v_1+iv_2=2i\varphi_{\bar{z}}, \quad\quad v_- = v_1-iv_2=-2i\varphi_z$$
for incompressible flow, preserving anti-holomorphicity constraint, we have dependence $\zeta=\zeta(\bar{z}, t)$ and the model reduces to the system
$$i\zeta_t -2\varphi_{\bar{z}}\zeta_{\bar{z}}+2\zeta_{\bar{z}\bar{z}}-\frac{4\bar{\zeta}\zeta^2_{\bar{z}}}{1+|\zeta|^2}=0,$$

\begin{equation}\label{DE5}
\varphi_{z\bar{z}}=-\frac{2\bar{\zeta}_z\zeta_{\bar{z}}}{(1+|\zeta|^2)^2}.
\end{equation}
We can rearrange the first equation as
$$i\zeta_t +2\zeta_{\bar{z}}\{-\varphi-2\ln(1+|\zeta|^2) + \ln\zeta_{\bar{z}}\}_{\bar{z}}=0$$
If we choose $\varphi=-2\ln(1+|\zeta|^2)$, then it is easy to verify the equation (\ref{DE5}) holds true automatically. Hence, we have the anti-holomorphic Schr\"odinger equation
$$i\zeta_t + 2\zeta_{\bar{z}\bar{z}}=0.$$
By the transformation $u=4(\ln\zeta)_{\bar{z}}$, which is a complex analog of the Cole-Hopf transformation, it implies the following complex Burgers' equation
$$iu_t+uu_{\bar{z}} + 2u_{\bar{z}\bar{z}}=0.$$
Here, $u$ can be interpreted as the complex velocity of effective flow with the complex potential $f(z) = \ln\zeta^4$. Then every zero of function $\zeta$ corresponds to the vortex solution (pole of complex velocity) of anti-holomorphic Burgers' equation. For more details we refer to \cite{GP, MPS} and reference therein.

\subsection{Vortex phenomena in higher dimensional cases}

Let us consider the case of Landau-Lifshitz equation from $\mathbb{R}^N$
with external magnetic fields (Schr\"odinger map flow equation) in the form
\begin{equation}
\label{L1.1}
\frac{\partial m}{\partial t}= m \times (\Delta m - m_3 \vec{e}_3) \quad \mbox{in} \ \R\times\R^N ,
\end{equation}
or equivalently the equation
\begin{equation}
\label{L1.2} -m \times \frac{\partial m}{\partial t}= \Delta m - m_3 \vec{e}_3+ (|\nabla m|^2 + m_3^2) m.
\end{equation}
Here $m : \R\times\R^N\to {\mathbb S}^2\subset\R^3$ and where $\vec{e}_3= (0,0, 1) \in \R^3$.
From the physical side of view, one expects topological solitons,
which are half magnetic bubbles, exist in solutions of (\ref{L1.1}) (see \cite{Hubert} and \cite{PS}).
Indeed, in \cite{HL1}-\cite{HL2}, F. Hang and F. Lin have established the corresponding static theory
for such magnetic vortices.

\medskip
By a reduction method, F. Lin and J. Wei \cite{linwei}
looked for a solution of the traveling wave $m(s', s_N-C t)$
(i.e., travel in the $s_N-$direction with the speed $C=c\epsilon >0$)
of the equation (\ref{L1.1}).
Then $m$ must be a solution of
\begin{equation}
\label{L1.3}
-c\epsilon \frac{\partial m}{\partial s_N}= m \times (\Delta m - m_3 \vec{e}_3).
\end{equation}
After a proper scaling in the space, (\ref{L1.3}) becomes
\begin{equation}
\label{L1.4}
-c\frac{\partial m}{\partial s_N}= m \times (\Delta m- \frac{m_3}{\epsilon^2} \vec{e}_3),
\quad
s \in \R^N,
\end{equation}
or
\begin{equation}
\label{L1.5}
c\,m \times \frac{\partial m}{\partial s_N }
= \Delta m -\frac{m_3 \vec{e}_3}{\epsilon^2} + ( |\nabla m|^2 +\frac{m_3^2}{\epsilon^2}) m.
\end{equation}
The main results of the paper \cite{linwei} are the following:
Let $N \geq 2$ and $\epsilon$ sufficiently small
there is an axially symmetric  solution $m=m(|s'|, s_N) \in C^\infty (\R^N, {\mathbb S}^2)$
of (\ref{L1.4}) such that
\begin{equation}
E_\epsilon (m)= \int_{\R^N} \frac{1}{2} (|\nabla m|^2 +\frac{m^2_3}{\epsilon^2} )\,{\mathrm d}s <\infty
\end{equation}
and that m has exactly one vortex at $ (|s'|, s_N)=(a_\epsilon, 0)$ of degree $+1$,
where $a_\epsilon \approx \frac{1}{2} $.
If $N=2$, the traveling velocity $C\sim \epsilon$, while $ C= (N-2) \epsilon | \log \frac{1}{\epsilon}|$ for $N\geq 3$.
Naturally such solution $m$ gives rise to a nontrivial (two-dimensional) traveling wave solution of (\ref{L1.1}) with a pair of vortex and antivortex which undergoes the Kelvin Motion as described  above.
Solutions constructed of this form are called {\em traveling vortex rings}
for the case of the dimension $N\geq 3$.
Later on, J. Wei and J. Yang \cite{weiyang} concerned the existence of traveling wave solutions
possessing vortex structures with rotating invariance for problem (\ref{L1.3}).

\subsection{Main results: Special solutions with vortices for (\ref{Ishimoriproblem}) and (\ref{Wavemaps})}
The problems on the existence and the blow-up profiles of wave maps are well-known
in analysis and attract much concerning in last few decades.
In particular, I. Rodnianski and his coauthor in \cite{RS, RSt} exhibit stable finite time blow up
regimes for the energy critical co-rotational wave map with the ${\mathbb S}^2$ target in all homotopy
classes and for the critical equivariant $SO(4)$ Yang-Mills problem.
They derive sharp asymptotics on the dynamics at blow up time and prove quantization of the energy focused at
the singularity. For the case of degree one, a sort of different blow-up profiles was shown in \cite{KST}.
The existence and regularity theory of wave maps from $\mathbb{R}^{1,2}$ was also established,
for more details we refer to \cite{ST} and reference therein.
However, it seems that one has not known whether or not (\ref{Wavemaps})
admits any solution with vortex structures by the knowledge of authors.

\medskip

Since we intend to seek for some special solutions of soliton type with vortex structure to (\ref{Ishimoriproblem})
and (\ref{Wavemaps}), let's recall that solitons are some special exact solutions
to classical non-linear field equations in physics. They are localized and have a finite energy.
In this sense, they behave like ordinary particles. Solitons and multi-solitons are stable because
they carry a topological charge $N$, which is an integer and equals the net number of particles.
Therefore, such solitons are usually called as topological solitons from the viewpoint of topology.
For instance, the Ishimori Model for the $su(2$) and $su(1, 1)$ algebras can be related to the ($2+1$)-dimensional
Davey-Stewartson equation with pseudo-euclidean space metric \cite{LMS, LS}.
However, in this representation nontrivial mappings from space-time to inner space,
classified by a topological charge, could arise.

\medskip
Because the topological charge is a conserved quantity, a single soliton cannot decay.
It is worthy to point out that the conservation of $N$ is not due to a Noether theorem,
but to the topological structure of the soliton.
On the other hand, some geometric Hamiltonian or dispersive flows are of subtle symmetry
when the starting manifolds and target manifolds are of invariance with respect to the groups of transformations.
Usually these geometric flows admit some symmetric solutions with respect to these Lie groups of transformations,
which obey some conservation laws due to Noether theorems and are of the properties of solitons.\medskip

We recall that in \cite{DingYin} Ding and Yin proposed to study the periodic solutions of the Schr\"odinger flow
in the case where the target manifold $N$ is a K\"ahler-Einstein manifold with positive scalar curvature.
such a class of special solutions can be viewed as a sort of geometric solitary wave solution (\cite{SongWang}).
Recently, from the viewpoint of differential geometry C. Song, X. Sun and Y. Wang \cite{SongSunWang}
propose a notion ``geometric solitons", which is concerning a special solutions for some geometric flows
and can be regarded as a geometric generalization of the classical solitary wave solutions.
More precisely,

\begin{definition}\cite{SongSunWang}
Suppose $(M, g)$ is a Riemannian manifold and $(N, h, J)$ is a K$\ddot{a}$hler manifold.
A solution u to the Schr\"odinger flow
\begin{align}\label{flow}
u_t = J(u)\tau(u).
\end{align}
is called a geometric soliton, if there exist two
one-parameter groups of isometries $\psi_t : M\rightarrow M$, $\phi_t : N\rightarrow N$
and a map $v : M\rightarrow N$ independent of $t$ such that $u$ has the form:
$$
u(t,x)\,=\,\phi_t\circ v\circ \psi_t.
$$
\end{definition}
For more explanations and examples for geometric solitons, the reader can refer to \cite{SongSunWang}.
In the present paper, we will use similar ideas to transfer problem (\ref{Ishimoriproblem})
and problem (\ref{Wavemaps}) to elliptic problems and then try to construct solutions with vortices.
For more details, we refer the reader to section \ref{section2}.

\medskip
Let ${\mathcal P}$ be the inverse of the stereographic projection
$$
\psi\,=\,\frac{ m_1+ i m_2}{1+m_3}:{\mathbb S}^2\rightarrow{\mathbb C},
$$
and then set
$$
{\mathcal R}(\tau, \cdot)={\mathcal P}\,e^{i\tau}.
$$
Here are main results.

\begin{theorem}\label{theorem-vortexpairs}
{\bf [Vortex Pairs]}
Assume that $N=2$.
\begin{enumerate}
  \item[(1).] For problem (\ref{Ishimoriproblem}), there is a solution in the form
$$
m(t,\tau, s_1,s_2)\,=\,{\mathcal R}\circ \psi(t,\tau, s_1,s_2)\in C^\infty (\R\times\R\times\R^2, {\mathbb S}^2),
$$
where $\psi(t, \tau,s_1,s_2)=U(s_1,s_2-c\tau-\omega t):\R\times\R^2\rightarrow{\mathbb C}$
possesses two vortices of degree $\pm 1$ traveling along the curves
$$
\{\,(x_1,x_2)\in \R^2\,:\, x_1=\pm d,\,  x_2\in\R\,\},
$$
under Kelvin motion with relation
$$
\frac{1}{d}\sim \frac{2c}{\sqrt{1-c^2}}-\frac{2\omega}{\sqrt{1-c^2}}.
$$
  \item[(2).] For problem (\ref{Wavemaps}), there is a solution in the form
$$
m(\tau, s_1,s_2)\,=\,{\mathcal R}\circ \psi(\tau, s_1,s_2)\in C^\infty (\R\times\R^2, {\mathbb S}^2),
$$
where $\psi(\tau,s_1,s_2)=U(s_1,s_2-c\tau):\R\times\R^2\rightarrow{\mathbb C}$
possesses two vortices of degree $\pm 1$ traveling along the curves
$$
\{\,(x_1,x_2)\in \R^2\,:\, x_1=\pm d,\,  x_2\in\R\,\},
$$
under Kelvin motion with relation
$$
\frac{1}{d}\sim \frac{2c}{\sqrt{1-c^2}}.
$$
\end{enumerate}
\qed
\end{theorem}

\begin{theorem}
\label{theorem-vortexrings}{\bf [Vortex rings]}
Assume that $N=3$.
\begin{enumerate}
\item[(1).] There is a solution to (\ref{Ishimoriproblem}) in the form
$$
m(t, \tau, s_1, s_2, s_3)\,=\,{\mathcal R}\circ \psi(t, \tau, s_1, s_2, s_3)\in C^\infty (\R\times\R\times\R^3, {\mathbb S}^2),
$$
where $\psi(t, \tau, s_1, s_2, s_3)=U(s_1, s_2, s_3-c\tau-\omega t):\R\times\R\times\R^3\rightarrow{\mathbb C}$
has a vortex helix directed along
the curve in the form
$$
\alpha\in\R\mapsto(d\cos\alpha, d\sin\alpha, 0)\in\R^3,
$$
and $u$ is also invariant under the skew motion expressed by cylinder coordinate
\begin{align*}
\Sigma:\,(r,\theta,s_3)\mapsto (r,\theta+\alpha, s_3),
\quad \forall\,\alpha\in\R.
\end{align*}
The traveling velocity and the geometric parameters of the vortex ring obeys the relation
$$
\frac{1}{d}\log d\sim \frac{4c}{\sqrt{1-c^2}}-\frac{4\omega}{\sqrt{1-c^2}}.
$$
\item[(2).] There is a solution to (\ref{Wavemaps}) in the form
$$
m(\tau, s_1,s_2)\,=\,{\mathcal R}\circ \psi(\tau, s_1,s_2,s_3)\in C^\infty (\R\times\R^3, {\mathbb S}^2),
$$
where $\psi(\tau,s_1,s_2,s_3)=U(s_1,s_2,s_3-c\tau):\R\times\R^3\rightarrow{\mathbb C}$
has a vortex helix directed along
the curve in the form
$$
\alpha\in\R\mapsto(d\cos\alpha, d\sin\alpha, 0)\in\R^3,
$$
and $u$ is also invariant under the skew motion expressed by cylinder coordinate
\begin{align*}
\Sigma:\,(r,\theta,s_3)\mapsto (r,\theta+\alpha, s_3),
\quad \forall\,\alpha\in\R.
\end{align*}
The traveling velocity and the geometric parameters of the vortex ring obeys the relation
$$
\frac{1}{d}\log d\sim \frac{4c}{\sqrt{1-c^2}}.
$$
\end{enumerate}
\qed
\end{theorem}

\medskip
\begin{remark}
Some words are in order to explain the results in the above theorems.
\begin{description}
  \item[(1).] For more information of the solutions, the reader can refer to Section \ref{section2}.

  \item[(2).] The reader can refer to Remarks \ref{remarkLocationandveloPairs} and \ref{remarkLocationandveloHelix}
for the relation between traveling velocity and locations of vortices.
In the proof these theorems, we need the traveling velocity is small enough
\end{description}
\end{remark}

\medskip
The construction of other vortex structures such as vortex helices is under working.
The organization of the paper is as follows: In subsection \ref{section21}, we introduce some transformations and then write problem (\ref{Wavemaps})
and problem (\ref{Ishimoriproblem}) into elliptic cases
for new unknown functions.
To describe the vortex phenomena (vortex pairs and vortex rings),
we then further formulate these problems by the symmetries in subsection \ref{section22}.
An outline of the strategies of proof will be provided in subsection \ref{section23}.
Some preliminaries are prepared in Section \ref{3preliminaries}.
Section \ref{section4} (Section \ref{section5}) is devoted to the construction
of solutions with vortex pairs (vortex rings), which will complete the proof of
Theorem \ref{theorem-vortexpairs} (Theorem \ref{theorem-vortexrings}).

\section{Various vortex structures: outline of the proofs}\label{section2}
\setcounter{equation}{0}
For the construction of various vortex structures and prove all theorems in Section \ref{section1},
the main objective of this section is to formulate the problems in suitable forms,
and then introduce some notation and describe the symmetries of the problems.
For the convenience of readers, we also provide an outline of the proofs in subsection \ref{section23}.

\subsection{Formulations of the problems}\label{section21}
\subsubsection{Wave maps}
\medskip
Setting the stereographic projection
$$
\psi\,=\,\frac{ m_1+ i m_2}{1+m_3},
$$
problem (\ref{Wavemaps}) becomes
\begin{equation}\label{steorographicproblemWave}
\Box \psi
\,-\,\frac{2 \bar{\psi}}{1+|\psi|^2}
\Big(\nabla_\tau\psi\cdot\nabla_\tau\psi-\nabla_s\psi \cdot \nabla_s\psi\Big)
\,=\,0.
\end{equation}
Here and throughout the paper, we use $ \bar{\psi}$ to denote the conjugate of $\psi$.

\medskip
We look for the traveling solitary wave solutions in the form
\begin{align*}
\begin{aligned}
\psi(\tau,s_1,\,\cdots,\,s_N)\,=\,U({\tilde s}_1,\,\cdots,\,{\tilde s}_N)\,e^{i\tau},
\\
\\
{\tilde s}_1=s_1,\,\cdots,\,{\tilde s}_{N-1}=s_{N-1},\,{\tilde s}_N=s_N-c\tau.
\end{aligned}
\end{align*}
This gives that $U$ will satisfy
\begin{align*}
\begin{aligned}
&\,i\,\frac{ 1-|U|^2}{1+|U|^2}\,
\,2\,c\,\frac{\partial U}{\,\partial {\tilde s}_N\,}
\,-\,c^2\frac{\,\partial^2 U\,}{\partial {\tilde s}_N^2}
\,+\,\frac{2{\bar U}}{1+|U|^2}
\,c^2\frac{\,\partial U\,}{\partial {\tilde s}_N}\frac{\,\partial U\,}{\partial {\tilde s}_N}
\\
&\,+\,\Delta U
\,-\,\frac{2 \bar{U}}{1+|U|^2}\nabla U \cdot \nabla U
\,+\,\frac{ 1-|U|^2}{1+|U|^2} U
\,=\,0.
\end{aligned}
\end{align*}
Suitable rescaling
\begin{align}
\begin{aligned}
U({\tilde s}_1,\,\cdots,\,{\tilde s}_N)\,=\,u({\hat s}_1,\,\cdots,\,{\hat s}_{N}),\qquad
\\
\\
{\tilde s}_1={\hat s}_1,\,\cdots,\,{\tilde s}_{N-1}={\hat s}_{N-1},\,{\tilde s}_{N}={\hat s}_{N}\sqrt{1-c^2}.
\end{aligned}
\end{align}
will derive that
\begin{align}\label{TraSolproblem-wave}
\begin{aligned}
&\,i\,\frac{ 1-|u|^2}{1+|u|^2}\,
\,\frac{2\,c}{\,\sqrt{1-c^2}\,}\,\frac{\partial u}{\,\partial {\hat s}_N\,}
\,+\,\Delta u
\,-\,\frac{2 \bar{u}}{1+|u|^2}\nabla u \cdot \nabla u
\,+\,\frac{ 1-|u|^2}{1+|u|^2} u
\,=\,0,
\end{aligned}
\end{align}
where we have denoted
\begin{align*}
\triangle u\,=\,
\sum_{j=1}^N\frac{\partial^2 u}{\partial {\hat s}_j^2},
\qquad&
\nabla u\,=\,\Big(\frac{\partial u}{\partial {\hat s}_1},\cdots,\frac{\partial u}{\partial {\hat s}_N}\Big).
\end{align*}
In the sequel, we will set
\begin{align}\label{F}
F(u)=\frac{ 1-|u|^2}{1+|u|^2} u.
\end{align}

\subsubsection{Generalized Schr\"odinger maps}

\medskip
Setting the stereographic projection
$$
\psi\,=\,\frac{ m_1+ i m_2}{1+m_3},
$$
problem (\ref{Ishimoriproblem}) becomes
\begin{equation}\label{steorographicproblemIshimori}
i\,\partial_t\psi
\,+\,\Box \psi
\,-\,\frac{2 \bar{\psi}}{1+|\psi|^2}
\Big(\nabla_\tau\psi\cdot\nabla_\tau\psi-\nabla_s\psi \cdot \nabla_s\psi\Big)
\,=\,0.
\end{equation}
Here and throughout the paper, we use $ \bar{\psi}$ to denote the conjugate of $\psi$.

\medskip
We look for the traveling solitary wave solutions in the form
\begin{align*}
\begin{aligned}
\psi(t, \tau, s_1,\,\cdots,\,s_N)\,=\,U({\tilde s}_1,\,\cdots,\,{\tilde s}_{N})\,e^{i\tau},
\\
\\
{\tilde s}_1=s_1,\,\cdots,\,{\tilde s}_{N-1}=s_{N-1},\,{\tilde s}_{N}=s_{N}-c\tau-\omega t.
\end{aligned}
\end{align*}
This leads to the problem
\begin{align*}
\begin{aligned}
&\,-\,i\,\omega\,\frac{\partial U}{\,\partial {\tilde s}_N\,}
\,+\,i\,\frac{ 1-|U|^2}{1+|U|^2}\,2\,c\,\frac{\partial U}{\,\partial {\tilde s}_N\,}
\,-\,c^2\frac{\,\partial^2 U\,}{\partial {\tilde s}_N^2}
\,+\,\frac{2{\bar U}}{1+|U|^2}
\,c^2\frac{\,\partial U\,}{\partial {\tilde s}_N}\frac{\,\partial U\,}{\partial {\tilde s}_N}
\\
&\,+\,\Delta U
\,-\,\frac{2 \bar{U}}{1+|U|^2}\nabla U \cdot \nabla U
\,+\,\frac{ 1-|U|^2}{1+|U|^2} U
\,=\,0.
\end{aligned}
\end{align*}
Suitable rescaling
\begin{align}
\begin{aligned}
U({\tilde s}_1,\,\cdots,\,{\tilde s}_N)\,=\,u({\hat s}_1,\,\cdots,\,{\hat s}_{N}),\qquad
\\
{\tilde s}_1={\hat s}_1,\,\cdots,\,{\tilde s}_{N-1}={\hat s}_{N-1},\,{\tilde s}_{N}={\hat s}_{N}\sqrt{1-c^2}.
\end{aligned}
\end{align}
will give that
\begin{align}\label{TraSolproblem-Ishimori}
\begin{aligned}
&-i\,\frac{\omega}{\,\sqrt{1-c^2}\,}\,\frac{\partial u}{\,\partial {\hat s}_N\,}
\,+\,
i\,\frac{ 1-|u|^2}{1+|u|^2}\,
\frac{2\,c}{\,\sqrt{1-c^2}\,}\,\frac{\partial u}{\,\partial {\hat s}_N\,}
\\
&\,+\,\Delta u
\,-\,\frac{2 \bar{u}}{1+|u|^2}\nabla u \cdot \nabla u
\,+\,\frac{ 1-|u|^2}{1+|u|^2} u
\,=\,0.
\end{aligned}
\end{align}

\medskip
\subsection{Vortex phenomena}\label{section22}
In this subsection, for future constructions of solutions with various vortex structures,
we shall further transform problem (\ref{TraSolproblem-wave})
and problem (\ref{TraSolproblem-Ishimori}) to suitable forms.
Note that we will assume that the traveling velocity parameters
$c$ and $\omega$ are small positive numbers in the present paper.

\subsubsection{Vortex pairs in two dimensional spaces: $N=2$}

We will construct solutions with pair of vortices to problem (\ref{TraSolproblem-wave}) and
problem (\ref{TraSolproblem-Ishimori}) in the case of $N=2$.
We will set
\begin{align}\label{vePairs}
\varepsilon\,=\,{2c\,}/{\sqrt{1-c^2}\,}.
\end{align}
Problem (\ref{TraSolproblem-wave}) can be written as
\begin{align}\label{orignalproblem-Q2-vortexpairs}
\begin{aligned}
i\,\varepsilon\,\frac{ 1-|u|^2}{1+|u|^2}
\frac{\partial u}{\,\partial x_2\,}
\,+\,\Delta u
\,+\, \frac{ 1-|u|^2}{1+|u|^2} u
\,-\,\frac{2 \bar{u}}{1+|u|^2}\nabla u \cdot \nabla u
\,=\,0.
\end{aligned}
\end{align}
We are also concerning problem (\ref{TraSolproblem-Ishimori}) in the form
\begin{align}\label{orignalproblem-Q2T2-vortexpairs}
\begin{aligned}
-i\,\kappa\,\varepsilon\,
\frac{\partial u}{\,\partial x_2\,}
\,+\,
i\,\varepsilon\,\frac{ 1-|u|^2}{1+|u|^2}
\frac{\partial u}{\,\partial x_2\,}
\,+\,\Delta u
\,+\, \frac{ 1-|u|^2}{1+|u|^2} u
\,-\,\frac{2 \bar{u}}{1+|u|^2}\nabla u \cdot \nabla u
\,=\,0,
\end{aligned}
\end{align}
where we have denoted
\begin{align}\label{kvePairs}
\kappa\,\ve=\omega/\sqrt{1-c^2}.
\end{align}
We will construct solutions with vortex pairs(a vortex and its antipair)  to problem (\ref{orignalproblem-Q2-vortexpairs})
and problem (\ref{orignalproblem-Q2T2-vortexpairs}) in Section \ref{section4}.

\medskip
\noindent\textbf{Notation and symmetry properties of the operators: }
For later use, we will denote
$$
z=x_1+ix_2
$$
and then set the notation
\begin{align}\label{S0QT}
\begin{aligned}
{\mathbb S}_0 [u]\,:=\,
\Delta u- \frac{2 \bar{u}}{1+|u|^2} \nabla u\cdot \nabla u \,+\, \frac{1-|u|^2}{1+|u|^2} u,
\\
\\
{\mathbb Q}_2 [u]\,:=\,
\,i\,\frac{1-|u|^2}{1+|u|^2}\,\frac{\partial u}{\partial x_2},
\qquad
%
\qquad
{\mathbb T}_2 [u]\,:=\,
-i\,\frac{\partial u}{\partial x_2},
\\
\\
{\mathbb S}_1 [u]
\,:=\,{\mathbb S}_0
\,+\,\ve\,{\mathbb Q}_2,
\quad
{\mathbb S}_2 [u]
\,:=\,{\mathbb S}_0
\,+\,\ve\,{\mathbb Q}_2
\,+\,\kappa\,\ve\,{\mathbb T}_2,
\end{aligned}
\end{align}
It is easy to see the operators
${\mathbb S}_1$ and ${\mathbb S}_2$
are invariant under the following two transformations
\begin{equation*}
u(z) \to \overline{ u(\bar{z})},
\quad
u (z) \to u(-\bar{z} ).
\end{equation*}
Thus we impose the following symmetry on the solutions $u$ to
problem (\ref{orignalproblem-Q2-vortexpairs}) and problem (\ref{orignalproblem-Q2T2-vortexpairs})
\begin{equation}
\label{symmetryvorxpairs}
\Pi\,:=\,\big\{ \,u(z)= \overline{u (\bar{z} )}, \ u(z)= u(-\bar{z}) \,\big\}.
\end{equation}
This symmetry will play an important role in our analysis.

\medskip
\subsubsection{Vortex Rings}

Here we only consider the case $N=3$ and then construct solutions
with vortex rings for problem (\ref{TraSolproblem-wave}) and problem (\ref{TraSolproblem-Ishimori}).
We set
\begin{align}\label{veHelix}
\varepsilon|\log\ve|=\frac{2\,c}{\,\sqrt{1-c^2}\,}.
\end{align}
By using the cylinder coordinates
$$
{\hat s}_1\,=\,r\cos\Psi,
\quad
{\hat s}_2\,=\,r\sin\Psi,
\quad
{\hat s}_3\,=\,{\hat s}_3,
$$
equation (\ref{TraSolproblem-wave}) is transferred to
\begin{align}
\begin{aligned}\label{inter2}
&i\, \ve\, |\log \ve |\, \frac{ 1-|u|^2}{1+|u|^2}\, \frac{\partial u}{\partial {\hat s}_3 }
\,+\,\frac{\partial^2 u}{\partial r^2}
\,+\,\frac{1}{r}\frac{\partial u}{\partial r}
\,+\,\frac{1}{r^2}\frac{\partial^2 u}{\partial \Psi^2}
\,+\,\frac{\partial^2 u}{\partial {\hat s}_3 ^2}
\,+\, \frac{ 1-|u|^2}{1+|u|^2} u
\\
&\quad\,=\,  \frac{2 \bar{u}}{1+|u|^2} \Bigg[\,\Big(\frac{\partial u}{\partial r}\Big)^2
+\frac{1}{r^2}\Big(\frac{\partial u}{\partial \Psi}\Big)^2
+\Big(\frac{\partial u}{\partial {\hat s}_3 }\Big)^2\Bigg].
\end{aligned}
\end{align}
For problem (\ref{inter2}), we want to find a solution $u$ which has a vortex ring directed along
the curve in the form
$$
\alpha\in\R\mapsto(d\cos\alpha, d\sin\alpha, 0)\in\R^3,
$$
with two parameters
\begin{align}\label{parameters}
d={\hat d}/\ve.
\end{align}
Moreover, $u$ is also invariant under the skew motion
\begin{align}
\Sigma:\,(r,\Psi,{\hat s}_3 )\mapsto (r,\Psi+\alpha, {\hat s}_3),\quad \forall\,\alpha\in\R.
\end{align}
Hence, $u$  satisfies the problem
\begin{align}
\begin{aligned}\label{inter3}
&i\,\ve\, |\log \ve |\frac{ 1-|u|^2}{1+|u|^2}  \frac{\partial u}{\partial {\hat s}_3 }
\,+\,\frac{\partial^2 u}{\partial r^2}
\,+\,\frac{1}{r}\frac{\partial u}{\partial r}
\,+\, \frac{\partial^2 u}{\partial {\hat s}_3 ^2}
\,+\, \frac{ 1-|u|^2}{1+|u|^2} u
\\
&\,=\,  \frac{2 \bar{u}}{1+|u|^2}\Bigg[\,\Big(\frac{\partial u}{\partial r}\Big)^2
+ \Big(\frac{\partial u}{\partial {\hat s}_3 }\Big)^2\Bigg].
\end{aligned}
\end{align}
For further convenience of notation, we also introduce the rescaling
\begin{equation}\label{coordinatesrs}
 (r, {\hat s}_3 )= ( x_1, x_2), \quad  z= x_1+ i x_2.
\end{equation}
Thus (\ref{inter3}) becomes
\begin{align}
\begin{aligned}
\label{orignalproblem-Q2-vortexhelices}
&i\ve |\log \ve | \frac{ 1-|u|^2}{1+|u|^2}\frac{\partial u}{\partial x_2}
\,+\,\frac{\partial^2 u}{\partial x_1^2}
\,+\,\frac{1}{x_1}\frac{\partial u}{\partial x_1}
\,+\, \frac{\partial^2 u}{\partial x_2^2}
\,+\, \frac{ 1-|u|^2}{1+|u|^2} u
\\
&\quad\,=\,  \frac{2 \bar{u}}{1+|u|^2}\Bigg[\,\Big(\frac{\partial u}{\partial x_1}\Big)^2
+ \Big(\frac{\partial u}{\partial x_2}\Big)^2\Bigg].
\end{aligned}
\end{align}

\medskip
We will also consider problem (\ref{TraSolproblem-Ishimori}) in the form
\begin{align}
\begin{aligned}\label{inter4}
-i\,\kappa\,\varepsilon\,|\log\ve|\,
\frac{\partial u}{\,\partial {\hat s}_3\,}
\,+\,
i\,\varepsilon\,|\log\ve|\,\frac{ 1-|u|^2}{1+|u|^2}
\frac{\partial u}{\,\partial {\hat s}_3\,}
\,+\,\Delta u
\,+\, \frac{ 1-|u|^2}{1+|u|^2} u
\,=\,
\frac{2 \bar{u}}{1+|u|^2}\nabla u \cdot \nabla u,
\end{aligned}
\end{align}
where
\begin{align}\label{kveHelix}
\kappa\,\varepsilon\,|\log\ve|=\omega/\sqrt{1-c^2}.
\end{align}
As we have done in the above, for the existence of solutions with vortex rings,
problem (\ref{inter4}) can be transformed to
\begin{align}\label{orignalproblem-Q2T2-vortexhelices}
\begin{aligned}
&-i\,\kappa\, \ve\, |\log \ve |\, \,\frac{\partial u}{\partial x_2}
\,+\,i \delta \ve |\log \ve | \frac{ 1-|u|^2}{1+|u|^2}\frac{\partial u}{\partial x_2}
\,+\,\frac{\partial^2 u}{\partial x_1^2}
\,+\,\frac{1}{x_1}\frac{\partial u}{\partial x_1}
\\
&\quad
\,+\, \frac{\partial^2 u}{\partial x_2^2}
\,+\, \frac{ 1-|u|^2}{1+|u|^2} u
\,=\,
\frac{2 \bar{u}}{1+|u|^2}\Bigg[\,\Big(\frac{\partial u}{\partial x_1}\Big)^2
+ \Big(\frac{\partial u}{\partial x_2}\Big)^2\Bigg].
\end{aligned}
\end{align}

\medskip
\textbf{The symmetric properties of the problems: }
\medskip
Some words are in order to describe the methods in dealing with
the existence of solutions of the above problems.
Problem
(\ref{orignalproblem-Q2-vortexhelices}),
and problem (\ref{orignalproblem-Q2T2-vortexhelices})
become two-dimensional cases.
The key point is then to construct a solution with vortices of degree $+1$ and their antipairs of degree $-1$.
Additional to the computations for standard vortices in two dimensional case,
there an  extra derivative term
\begin{align}\label{singularity0}
\frac{1}{x_1} \frac{\partial u}{\partial x_1}.
\end{align}
In the following we shall mainly focus on how to deal with these terms by the methods in \cite{linwei}.

\medskip
By using the symmetries, in the sequel, we shall consider problems
(\ref{orignalproblem-Q2-vortexhelices}),
(\ref{orignalproblem-Q2T2-vortexhelices})
with imposing the boundary conditions
\begin{align}\label{boundary}
\begin{aligned}
|u(z)| \to 1 \quad \mbox{as} \ |z| \to +\infty,\qquad\quad
\\
\frac{\partial u}{\partial x_1} (0, x_2)\,=\,0,\quad
\forall \,x_2\in\R.
\end{aligned}
\end{align}
Moreover, it is easy to see that problem (\ref{orignalproblem-Q2-vortexhelices})
and problem (\ref{orignalproblem-Q2T2-vortexhelices}) are also invariant under the following two transformations
\begin{equation}
 u(z) \to \overline{ u(\bar{z})}, \ \ u (z) \to u(-\bar{z} ).
\end{equation}
Thus we impose the following symmetry on the solution $u$
\begin{equation}
\label{symmetryvortexhelix}
\Pi :=\Big\{\, u(z)= \overline{u (\bar{z} )}, \ u(z)= u(-\bar{z}) \,\Big\}.
\end{equation}
This symmetry will play an important role in our analysis.
As a conclusion, if we write the solutions $u$ to
problem (\ref{orignalproblem-Q2-vortexhelices}) and problem (\ref{orignalproblem-Q2T2-vortexhelices})
in the form
$$
u(x_1,x_2)=u_1(x_1,x_2)+iu_2(x_1,x_2),
$$
then $u_1$ and $u_2$ enjoy the following conditions:
\begin{align}
\begin{aligned}\label{symmetryandboundary}
u_1(x_1,x_2)=u_1(-x_1,x_2),
&\qquad
u_1(x_1,x_2)=u_1(x_1,-x_2),
\\
u_2(x_1,x_2)=u_2(-x_1,x_2),
&\qquad
u_2(x_1,x_2)=-u_2(x_1,-x_2),
\\
\frac{\partial u_1}{\partial x_1}(0,x_2)=0,
&\qquad
\frac{\partial u_2}{\partial x_1}(0,x_2)=0.
\end{aligned}
\end{align}
Note that we shall construct the approximate solution
to satisfy the boundary conditions in (\ref{symmetryandboundary}).
We will use these conditions to construct solutions to problem (\ref{orignalproblem-Q2-vortexhelices})
and problem (\ref{orignalproblem-Q2T2-vortexhelices}) in Section \ref{section5}.

\medskip
Recall the notation in (\ref{S0QT}).
We denote
\begin{align}
\begin{aligned}\label{H-operators}
{\mathbb H}_1 [u]:\,=\,
\frac{1}{x_1} \frac{\partial u}{\partial x_1}.
\end{aligned}
\end{align}
For later use, we also set
\begin{align}
\begin{aligned}\label{S3-S4}
{\mathbb S}_3 [u]
\,=\,{\mathbb S}_0
\,+\,\ve\,|\log\ve|\,{\mathbb Q}_2
\,+\,{\mathbb H}_1,\qquad
\\
\\
{\mathbb S}_4 [u]
\,=\,{\mathbb S}_0
\,+\,\ve\,|\log\ve|\,{\mathbb Q}_2
\,+\,\kappa\,\ve\,|\log\ve|\,{\mathbb T}_2
\,+\,{\mathbb H}_1.
\end{aligned}
\end{align}

\vspace{0.1cm}
\subsection{Outline of the Proof}\label{section23}
To prove all theorems in Section \ref{section1},
we will use the finite dimensional reduction
method to find solutions to
problems (\ref{orignalproblem-Q2-vortexpairs}),
(\ref{orignalproblem-Q2T2-vortexpairs}) in Section \ref{section4},
and also for problems (\ref{orignalproblem-Q2-vortexhelices}), (\ref{orignalproblem-Q2T2-vortexhelices})
in Section \ref{section5}.
The finite dimensional reduction procedure has been used in many other problems.
See \cite{DFM}, \cite{GW}, \cite{LNW} and the references therein.
M. del Pino, M. Kowalczyk and M. Musso \cite{DKM}  were the first
to use this procedure to study Ginzburg-Landau equation in a bounded domain.
The methods in the present paper basically follow those in dealing with vortex phenomena
for Shr$\ddot{o}$dinger map equations in \cite{linwei}.
Here are the main steps of the approach.

\medskip
\noindent
{\bf Step 1: Constructions of approximate solutions}

\medskip
To construct a real solution,
the first step is to construct an approximate solution, denoted by $V_d$ in (\ref{Vdpairs}) or (\ref{Vd-helix}),
possessing a vortex located at $\vec{e}_1=(d,0)$
and its antipair  located at $\vec{e}_2=(-d,0)$.
Here $d$ is a parameter to be determined in the reduction procedure.
In the construction of approximation to solutions to problem
(\ref{orignalproblem-Q2-vortexhelices})
and problem (\ref{orignalproblem-Q2T2-vortexhelices}),
there are singularities caused by the application of the terms in (\ref{singularity0})
to the phase terms of the standard vortices,
which will be described in subsection \ref{subsection51}.

\medskip
The approximate solution $V_d$ has the symmetry
\begin{eqnarray}\label{symmetry0}
u_2(x_1,x_2)=\overline{u_2(x_1,-x_2)},
\qquad
u_2(x_1,x_2)=u_2(-x_1,x_2).
\end{eqnarray}
By substituting $V_d$ into problems
(\ref{orignalproblem-Q2-vortexpairs}),
(\ref{orignalproblem-Q2T2-vortexpairs}),
(\ref{orignalproblem-Q2-vortexhelices}),
(\ref{orignalproblem-Q2T2-vortexhelices}),
we can derive the estimations of
the errors such as ${\mathbb S}_1[V_d]$, ${\mathbb S}_2[V_d]$, ${\mathbb S}_4[V_d]$   in suitable weighted norms.
The reader can refer to the papers \cite{DKM} and \cite{linwei}.

\medskip
\noindent
{\bf Step 2: Finding a perturbation }

\medskip
We intend to look for solutions to problems
(\ref{orignalproblem-Q2-vortexpairs}),
(\ref{orignalproblem-Q2T2-vortexpairs}),
(\ref{orignalproblem-Q2-vortexhelices}),
(\ref{orignalproblem-Q2T2-vortexhelices}),
by adding a perturbation term, say $\psi$,  to the approximation $V_d$
where the perturbation term is small in suitable norms.
More precisely, for the perturbation $\psi=\psi_1+i\psi_2$ with symmetries
in (\ref{symmetryvorxpairs}) or (\ref{symmetryvortexhelix}) ,
we take the solution $u$ in the form
(cf. (\ref{perturbationdecompositionpairs}) or (\ref{perturbationdecompositionhelix}))
$$
 u (y)= \eta ( V_d + i V_d \psi) + (1-\eta) V_d e^{i \psi}.
$$
This perturbation method near the vortices was introduced in \cite{DKM}.

\medskip
For given parameters $d$ and $\ve$ small,
instead of considering the problems
(\ref{orignalproblem-Q2-vortexpairs}),
(\ref{orignalproblem-Q2T2-vortexpairs}),
(\ref{orignalproblem-Q2-vortexhelices}),
and (\ref{orignalproblem-Q2T2-vortexhelices})
we look for a $\psi$ to their projected forms
in (\ref{nonlinearprojection-pairs}) or (\ref{nonlinearprojection-helix}).
By writing the projected problem
in the form of the perturbation term $\psi$(with a linear part and a nonlinear part),
we can find the perturbation term $\psi$ through {\em a priori estimates} and the contraction mapping theorem.

\medskip

\noindent
{\bf Step 3: Adjusting the parameters }

\medskip
Note that the perturbation term $\psi$ and the Lagrange multipliers $c$ in (\ref{nonlinearprojection-pairs})
or (\ref{nonlinearprojection-helix}) are functions
of the parameters $d$.
To get real solutions to problems
(\ref{orignalproblem-Q2-vortexpairs}),
(\ref{orignalproblem-Q2T2-vortexpairs}),
(\ref{orignalproblem-Q2-vortexhelices}),
and (\ref{orignalproblem-Q2T2-vortexhelices}),
we shall choose suitable parameters $d$ such that $c$ is zero.
It is equivalent to solve a reduced algebraic system for the Lagrange multiplier
\begin{eqnarray*}
\label{reduced0}
c_\ve (d)=0.
\end{eqnarray*}
We can derive the equations in (\ref{reduced0}) by the standard reduction procedure.
This will be done in subsections \ref{reductionprocedurepairs} or \ref{reductionprocedurehelix}.
In other words, we achieve the balance between
the vortex-antivortex interaction and the effect of motion of vortices by
adjusting the locations of the vortices.

\section{Some preliminaries}\label{3preliminaries}
\setcounter{equation}{0}

In this section, we collect some important facts which will be used later.
These include the asymptotic behaviors and nondegeneracy of degree one vortex.

\medskip
For the problem
\begin{align}\label{modelproblem}
\Delta u
\,+\, \frac{ 1-|u|^2}{1+|u|^2} u
\,-\,\frac{2 \bar{u}}{1+|u|^2}\nabla u \cdot \nabla u
\,=\,0
\quad\mbox{on}\quad\R^2,
\end{align}
if we look for  solutions of standard vortex of degree $+1$ in the form
$$
u= w^{+} :=\rho (\ell) e^{i \theta},
$$
in polar coordinate $(\ell,\theta)$,  then  $\rho$ satisfies
\begin{equation}
\label{8}
\rho''
\,+\,\frac{\rho'}{\ell}
\,-\,\frac{2 \rho (\rho')^2}{ 1+\rho^2}
\,+\, \Big(\,1-\frac{1}{\ell^2}\,\Big) \frac{1-\rho^2}{1+\rho^2} \rho\,=\,0.
\end{equation}
Another solution $ w^{-}:= \rho(\ell) e^{-i \theta}$ will be of vortex of degree $-1$.
These two functions will be our  block elements for future construction of approximate solutions.

\medskip
\noindent{\bf Notation: }
{\em For simplicity, from now on,  we use $
w=\rho(\ell)e^{i\theta} $ to denote the degree $+1$ vortex.
\qed}

\medskip
 The following properties of $\rho$ are proved in \cite{HL1}.
\begin{lemma}
There hold the asymptotic behaviors:
\\
${\mathrm (1)}\,$ $\rho(0)=0,\quad 0<\rho (\ell) <1,\quad \rho' >0$ for $ \ell>0$,
\\
${\mathrm (2)}\,$ $ \rho (\ell)=1 - c_0 \frac{e^{-\ell}}{\sqrt{\ell}} + O( \ell^{-3/2} e^{-\ell})$ as $ \ell \to +\infty$, where $c_0>0$.
\qed
\end{lemma}

\medskip
Setting $w=w_1+iw_2$ and $z=x_1+ix_2$,
then we need to study  the following linearized problem of (\ref{modelproblem})
around the standard profile $w$:
\begin{align*}
{\mathbb L}_0 (\phi) \,=\,&
\Delta\phi
\,-\, \F{4(w_1\nabla w_1+ w_2\nabla w_2)}{1+|w|^2}\nabla\phi
\,-\,\F{4\nabla\langle w,\phi\rangle}{1+|w|^2}\nabla w
\\
&\,+\,\F{8\langle w,\phi\rangle(w_1\nabla w_1+ w_2\nabla w_2)}{(1+|w|^2)^2}\nabla w
\,+\,\F{4\langle \nabla w,\nabla \phi\rangle}{1+|w|^2}w
\\
&\,-\, \F{4(1+|\nabla w|^2)\langle w,\phi\rangle}{(1+|w|^2)^2}w
\,+\,\F{2|\nabla w|^2}{1+|w|^2}\phi
\,+\,\F{1-|w|^2}{1+|w|^2}\phi.
\end{align*}
The nondegeneracy of $w$ is contained in the following lemma\cite{linwei}.

\begin{lemma}\label{non1}
Suppose that
\begin{equation}
\label{L1n}
{\mathbb L}_0 [\phi]=0, \ \
\end{equation}
where $\phi= iw \psi$, and $ \psi=\psi_1+i \psi_2$ satisfies the following decaying estimates
\begin{align}
\begin{aligned}
\label{psidecay}
|\psi_1|\,+\,|z||\nabla \psi_1|  &\,\leq\, C (1+|z|)^{-\varrho},
\\
|\psi_2|\,+\,|z| |\nabla \psi_2|  &\,\leq\, C (1+|z|)^{-1-\varrho},
\end{aligned}
\end{align}
for  some $ 0<\varrho<1$.
Then
\begin{align*}
\phi=c_1 \F{\partial w}{\partial x_1}+c_2 \F{\partial w}{\partial x_2}
\end{align*}
for certain real constants $c_1$, $c_2$. \qed
\end{lemma}

\medskip
Finally,  we  will need the following lemma on decay estimates of a linear problem in $\R^2$.
\begin{lemma}
\label{l5.2}
Let $h$ satisfy
\begin{equation}
\Delta h+ f(z)=0,  h(\bar{z})=-h(z), |h| \leq C
\end{equation}
where $f$ satisfies
\begin{equation}
\label{fdecay}
 |f(z) | \leq \frac{C}{ (1+|z|)^{2+\varrho}}, 0<\varrho<1.
\end{equation}
Then
\begin{equation}
|h(z)| \leq \frac{C}{(1+|z|)^{\varrho}}.
\end{equation}
\qed
\end{lemma}
For the proof of the last two lemmas, the reader can refer to \cite{linwei}.

\section{Vortex pairs: construction of solutions to
(\ref{orignalproblem-Q2-vortexpairs}) and (\ref{orignalproblem-Q2T2-vortexpairs}) }
\label{section4}
\setcounter{equation}{0}

In this section, we will construct solutions with vortex pairs
to problem (\ref{orignalproblem-Q2-vortexpairs}) and problem (\ref{orignalproblem-Q2T2-vortexpairs}),
which will complete the proof for Theorem \ref{theorem-vortexpairs}.

\subsection{Approximate Solutions}

Using the degree one vortex, we introduce  the {\em approximate solutions}.
Throughout the paper, we assume that the distance $d$ satisfies
\begin{equation}
d=\frac{\hat{d}}{\ve}, \ \ \ \mbox{where} \ \  \hat{d} \in \big[\,{1}/{100},\, 100\,\big].
\end{equation}
The approximate function is then defined by
\begin{equation}
\label{Vdpairs}
V_d (z):= w^{+} (z- \vec{e}_1) w^{-} (z-\vec{e}_2 ),
\end{equation}
where ${\vec e}_1=(d,0)$, ${\vec e}_2=(-d,0)$.
It is easy to see that $ V_d (z) \in \Pi$ with $\Pi$ defined in (\ref{symmetryvorxpairs}).
Furthermore, a simple computation shows that
\begin{equation}
|V_d(z)| \to 1 \ \mbox{as} \ |z| \to +\infty.
\end{equation}
In fact,  for $|z| >>d$, we have
\begin{eqnarray*}
V_d &\approx & e^{i \theta_{\vec{e}_1} - i \theta_{\vec{e}_2}}
\\
& \approx & \frac{ (x_1^2 -d^2 + x_2^2 +2i x_2 d)}{\, \sqrt{ (x_1-d)^2 +x_2^2}\, \sqrt{ (x_1+d)^2 +x_2^2} \,}
\\
& \approx & 1 + O\Big(\frac{d^2}{|z|^2} +\frac{d}{|z|} \Big).
\end{eqnarray*}
Here $\theta_{\xi}$ denotes the angle argument around $ \xi$ and $ \ell_\xi=|z-\xi|$.
It is easy to see that
\begin{equation}
\nabla \ell_{\xi}= \frac{1}{|z-\xi|} (z-\xi),
\quad
\nabla \theta_{\xi}= \frac{1}{|z-\xi|^2} (z-\xi)^\perp,
\quad
|\nabla \theta_\xi |=\frac{1}{|z- \xi|},
\end{equation}
where  we denote $ z^\perp= (-x_2, x_1)$.

\subsection{Error Estimates}

We plug in the approximate function (\ref{Vdpairs}) into the operators
$ {\mathbb S}_1$  and ${\mathbb S}_2$ and obtain the errors
\begin{align}
\begin{aligned}
{\mathbb E}_1\,:=\,&\, {\mathbb S}_1 [V_d]\,=\, {\mathbb S}_0[V_d] \,+\, \ve\,{\mathbb Q}_2[V_d],
\\
\\
{\mathbb E}_2\,:=\,&\, {\mathbb S}_2 [V_d]\,=\, {\mathbb S}_0[V_d]
\,+\, \ve\,{\mathbb Q}_2[V_d]
\,+\, \kappa\,\ve\,{\mathbb T}_2[V_d].
\end{aligned}
\end{align}
Our purpose in this subsection is to estimate these two errors.

\medskip
By our construction and the properties of $\rho$, we have
\begin{equation}
V_d= e^{i\theta_{\vec{e}_1} -i\theta_{\vec{e}_2}}
\Big[\,1+ O\big(e^{- \min  (|z-\vec{e}_1|,\, |z-\vec{e}_2|)}\big)\,\Big].
\end{equation}
We divide $\R^2$ into two regions:
$$
\R_{+}^2\,=\,\big\{ x_1>0\big\},
\quad
\R_{-}^2\,=\,\big\{ x_1<0\big\}.
$$
By the symmetry assumption in (\ref{symmetryvorxpairs}), we just need to consider the region $\R^2_{+}$.

\medskip
In the region $\R^2_{+}$, we consider two cases. First, if $ |z-\vec{e}_2| > \frac{d}{2}+ \frac{|z-\vec{e}_1|}{2}$, then we have
\begin{eqnarray*}
 V_d&=& w(z-\vec{e}_1) e^{-i \theta_{\vec{e}_2}} \Big[\,1+ O\big(e^{-d/2 -|z-d\vec{e}_1|/2 }\big)\,\Big],
\\
\\
 |V_d| &=&  \rho (|z-\vec{e}_1|) \Big[\,1+ O\big(e^{-d/2-|z-\vec{e}_1|/2}\big)\,\Big].
\end{eqnarray*}
Secondly, if $ |z-\vec{e}_2| < \frac{d}{2}+ \frac{|z-\vec{e}_1|}{2}$,  we use the translated variable as follows
\[ z\,=\, \vec{e}_1 + y.\]
Then we obtain the following estimates
\begin{eqnarray*}
\frac{1-|V_d|^2}{1+|V_d|^2} V_d &=& \frac{1-\rho^2}{1+\rho^2}\, w\, \Big[\,1+O\big(e^{-d/2 -|y|/2}\big)\,\Big]
e^{-i \theta_{\vec{e}_2}},
\\
\\
 \Delta V_d &=& \left[\, \Delta w -2 i \nabla w \cdot \nabla \theta_{\vec{e}_2}
 - \frac{w}{|y+ \vec{e}_1-\vec{e}_2|^2} + O\big(e^{-d/2-|y|/2}\big)\, \right] e^{-i \theta_{\vec{e}_2}},
\\
\\
 \frac{2 \bar{V_d}}{1+|V_d|^2} \nabla V_d \cdot \nabla V_d
 &= & \frac{2 \bar{w}}{1+\rho^2}\, e^{-i \theta_{\vec{e}_2}}
 \times \Bigg[ \nabla w \cdot \nabla w -2i w \nabla w \cdot \nabla \theta_{\vec{e}_2} - w^2 \nabla \theta_{\vec{e}_2} \cdot \nabla \theta_{\vec{e}_2}  \Bigg]
\\
\\
&&\,+\, O(e^{-d/2-|y|/2}).
\end{eqnarray*}
Combining the estimates above, we have  for $ z \in \R_{+}^2$,
\begin{equation}
\label{S0(Vd)}
 {\mathbb S}_0 [V_d]=  e^{-i \theta_{\vec{e}_2}} \Bigg[ \frac{2 (\rho^2-1)}{\rho^2+1} i \nabla w \cdot \nabla \theta_{\vec{e}_2}
  + \frac{\rho^2 -1}{\rho^2+1 } \frac{w}{\,|y+\vec{e}_1-\vec{e}_2|^2\,} + O(e^{-d/2-|y|/2}) \Bigg].
\end{equation}

\medskip
On the other hand, we can estimate the term ${\mathbb Q}_2 [V_d]$ as follows
\begin{align}
\begin{aligned}\label{Q2(Vd)}
{\mathbb Q}_2 [V_d]
\,=\,&\,
i\, \frac{\partial V_d}{\partial x_2}\,\frac{1-|V_d|^2}{1+|V_d|^2}
\\
\,=\, &\,i\,\frac{\partial}{\partial y_2} \Bigg[ w (y) e^{-i \theta_{\vec{e}_2}} + O( e^{-d/2-|y|/2})\Bigg]
\frac{1-\rho^2}{1+\rho^2}
\Big(1+O( e^{-d/2-|y|/2})\Big)
\\
\,=\,&\,
i\, \Bigg[ \frac{\partial w}{\partial y_2} -i w \frac{\partial \theta_{\vec{e}_2}}{\partial y_2} \Bigg]
\frac{1-\rho^2}{1+\rho^2}
\,e^{-i \theta_{\vec{e}_2}} + O(e^{-d/2-|y|/2} ).
\end{aligned}
\end{align}
Similarly,
\begin{align}
\begin{aligned}\label{T2(Vd)}
{\mathbb T}_2 [V_d]
\,=\,
-i\, \frac{\partial V_d}{\partial x_2}
\,=\,
-i\,\Bigg[ \frac{\partial w}{\partial y_2} -i w \frac{\partial \theta_{\vec{e}_2}}{\partial y_2} \Bigg]
\,e^{-i \theta_{\vec{e}_2}} + O( e^{-d/2-|y|/2} ).
\end{aligned}
\end{align}
Note that
\begin{align}
\begin{aligned}
\frac{\partial }{\partial y_2} \theta_{\vec{e}_2}
= \frac{ y_1+2d}{(y_1+2d)^2+y_2^2}= O\Big(\frac{1}{d}\Big)= O(\ve),
\\
\\
\frac{\partial^2 }{\partial y^2_2} \theta_{\vec{e}_2}
= O\Big(\frac{1}{d^2}\Big)= O(\ve^2).\qquad
\end{aligned}
\end{align}

\medskip
In summary, we have obtained for $ z \in \R_{+}^2,\, z\,=\,\vec{e}_1 +y$
\begin{align}\label{S1(Vd)}
\begin{aligned}
{\mathbb E}_1\,=\,{\mathbb S}_1 [V_d]&\,=\,{\mathbb S}_0[V_d]\,+\,\ve\,{\mathbb Q}_2[V_d]
\\
&\,=\,
\Bigg[ \frac{2 (\rho^2-1)}{\rho^2+1} i \nabla w \cdot \nabla \theta_{\vec{e}_2}
+ \frac{\rho^2 -1}{\rho^2+1 } \frac{w}{|y+\vec{e}_1-\vec{e}_2|^2}\Bigg]e^{-i \theta_{\vec{e}_2}}
\\
&\quad\,+\,i\,\ve\, \Bigg[ \frac{\partial w}{\partial y_2} -i w \frac{\partial \theta_{\vec{e}_2}}{\partial y_2} \Bigg]
\frac{1-\rho^2}{\,1+\rho^2\,}\,e^{-i \theta_{\vec{e}_2}},
\end{aligned}
\end{align}
\begin{align}\label{S2(Vd)}
\begin{aligned}
{\mathbb E}_2\,=\,{\mathbb S}_2 [V_d]&\,=\,{\mathbb S}_0[V_d]
\,+\,\ve\,{\mathbb Q}_2[V_d]
\,+\,\kappa\,\ve\,{\mathbb T}_2[V_d]
\\
&\,=\,
\Bigg[ \frac{2 (\rho^2-1)}{\rho^2+1} i \nabla w \cdot \nabla \theta_{\vec{e}_2}
+ \frac{\rho^2 -1}{\rho^2+1 } \frac{w}{|y+\vec{e}_1-\vec{e}_2|^2}\Bigg]e^{-i \theta_{\vec{e}_2}}
\\
&\quad\,+\,i\,\ve\, \Bigg[ \frac{\partial w}{\partial y_2} -i w \frac{\partial \theta_{\vec{e}_2}}{\partial y_2} \Bigg]
\frac{1-\rho^2}{\,1+\rho^2\,}\,e^{-i \theta_{\vec{e}_2}}
\,-\,i\,\kappa\,\ve\, \Bigg[ \frac{\partial w}{\partial y_2} -i w \frac{\partial \theta_{\vec{e}_2}}{\partial y_2} \Bigg]
\,e^{-i \theta_{\vec{e}_2}}.
\end{aligned}
\end{align}
Similar (and almost identical) estimates  also hold in the region $\R^2_{-}$.

\subsection{Setting-up of the Problem}

We look for solutions to (\ref{orignalproblem-Q2-vortexpairs})
and (\ref{orignalproblem-Q2T2-vortexpairs}) in the form
\begin{equation}
\label{perturbationdecompositionpairs}
 u (y)= \eta ( V_d + i V_d \psi) + (1-\eta) V_d e^{i \psi}
\end{equation}
where $\eta$ is a function such that
\begin{equation}
\label{etavortexpairs}
 \eta= \tilde{\eta} (|z-\vec{e}_1|)+ \tilde{\eta} (|z-\vec{e}_2|).
\end{equation}
In the above $\tilde\eta$ is the cut-off function defined by
\begin{align}\label{cutoff}
\tilde{\eta} (s)=1\quad\mbox{for}\quad s \leq 1,
\qquad
\tilde{\eta} (s)=0\quad\mbox{for}\quad s \geq 2.
\end{align}
We may write $\psi=\psi_1+ i \psi_2$ with $\psi_1, \psi_2$ real-valued.
The symmetry imposed on $u$ (see (\ref{symmetryvorxpairs})) can be transmitted to the symmetry on $\psi$
\begin{equation}
\label{psisymvirtexpairs}
\psi (\bar{z})=-\overline{\psi (z)}, \ \ \psi(z)= \psi (- \bar{z} ).
\end{equation}
This symmetry will be of  importance in solving the linear problems.

\medskip
By using (\ref{orignalproblem-Q2-vortexpairs}),
we can formulate the equation for the perturbation term $\psi$.
Let $R>1$ be fixed positive constant.
In the inner region
$$
{\mathfrak S}_1\,=\,\Big\{\, z\,|\,z \in B_{9R} (\vec{e}_1) \cup B_{9R} (\vec{e}_2)\,\Big\},
$$
by setting
\begin{equation}
\phi= \eta i V_d \psi + (1-\eta) V_d (e^{i \psi} -1),
\end{equation}
we have
\begin{equation}
u=V_d +\phi
\end{equation}
and then write the equation (\ref{orignalproblem-Q2-vortexpairs}) in $\phi$ as
\begin{equation}
{\mathbb L}_1 [\phi]
\,+\, {\mathbb N}_1 [\phi]
\,+\, {\mathbb M}_1 [\phi] \,=\, {\mathbb E}_d,
\end{equation}
where we have denoted
$$
{\mathbb E}_d=-{\mathbb S}_1 [V_d],
$$
\begin{align*}
{\mathbb L}_1[\phi]=\Delta \phi -\frac{4 \bar{V_d}}{1+|V_d|^2} \nabla V_d \cdot \nabla \phi
-\frac{2 \bar{\phi}}{1+|V_d|^2} \nabla V_d \cdot \nabla V_d + \frac{2 \bar{V_d} (V_d \bar{\phi}
+ \bar{V_d} \phi)}{(1+|V_d|^2)^2} \nabla V_d \cdot \nabla V_d  + F' (V_d) \phi.
\end{align*}
The nonlinear operator is
\begin{align}
\begin{aligned}\label{nonlinearoperator-vortexpairs}
{\mathbb N}_1[\phi]\,=\,&\, F(V_d+\phi)- F(V_d)
\,-\,F' (V_d) \phi
\,+\, O((1+|\phi|) |\nabla \phi|^2),
\\
\\
{\mathbb M}_1[\phi]\,=\,&\,
\ve\,{\mathbb Q}_2[V_d+\phi]
\,-\,\ve\,{\mathbb Q}_2[V_d].
\end{aligned}
\end{align}
Note that, for problem (\ref{orignalproblem-Q2T2-vortexpairs})
we have a similar formulation, but the nonlinear operator ${\mathbb M}_1$ is given by
\begin{align*}
{\mathbb M}_1[\phi]\,=\,&\,
\ve\,{\mathbb Q}_2[V_d+\phi]
\,-\,\ve\,{\mathbb Q}_2[V_d]
\,+\,\kappa\,\ve\,{\mathbb T}_2[V_d+\phi]
\,-\,\kappa\,\ve\,{\mathbb T}_2[V_d].
\end{align*}

\medskip
We now formulate the problem for the perturbation term $\psi$ in the outer region
$$
{\mathfrak S}_2\,=\,\Big\{\, z\,|\, z\in \big(B_{4R} (\vec{e}_1) \cup B_{4R} (\vec{e}_2)\big)^c\,\Big\},
$$
in which we have $ u= V_d e^{i \psi}$.
By simple computations we obtain
\begin{align*}
\frac{\, {\mathbb S}_1 [V_d e^{i \psi}]\,}{iV_d e^{i \psi}}
\,=\,&\,
\Delta \psi
\,+\, 2 \frac{\, 1-|V_d|^2 + |V_d|^2 (e^{-2\psi_2} -1)\,}{ V_d (1+|V_d|^2 e^{-2\psi_2})} \nabla V_d \cdot \nabla \psi
\,-\,i\, \Big(\frac{2 |V_d|^2}{1+|V_d|^2 e^{-2\psi_2}} -1\Big) \nabla \psi \cdot \nabla \psi
\\
&
\,+\, \frac{1}{\,i V_d\,} \frac{2 |V_d|^2 \bar{V_d} (e^{-2\psi_2} -1)}{\, (1+|V_d|^2)(1+|V_d|^2 e^{-2\psi_2})\,}
\,\nabla V_d \cdot \nabla V_d
\,-\,i\, \frac{ 2|V_d|^2 (1-e^{-2\psi_2})}{\, (1+|V_d|^2 e^{-2 \psi_2})(1+|V_d|^2)\,}
\\
&\,+\,2\,\ve\,\frac{\partial V_d}{\partial x_2}\frac{ 2|V_d|^2 (1-e^{-2\psi_2})}{\, (1+|V_d|^2 e^{-2 \psi_2})(1+|V_d|^2)\,V_d\,}
\,+\,2\,i\,\ve\,\frac{\partial \psi}{\partial x_2}\frac{ 1-2|V_d|^2e^{-2\psi_2}}{\, 1+|V_d|^2 e^{-2 \psi_2}\,\,}
\,+\, \frac{\,{\mathbb E}_1\,}{\,i V_d\,}\, .
\end{align*}
The above equation can be also formulated in the form
\begin{equation}
{\mathbb L}_2 [\psi] \,+\, {\mathbb N}_2 [\psi] \,+\,{\mathbb M}_2 [\psi]\,=\,{E},
\end{equation}
where we have denoted
\begin{eqnarray*}
{\mathbb L}_2 [\psi] &=&
\Delta \psi
\,+\, \frac{2 (1-|V_d|^2)}{1+|V_d|^2} \nabla V_d \cdot \nabla \psi
\,+\, \frac{ 4i \bar{V_d}^2 \psi_2}{(1+|V_d|^2)^2} \nabla V_d \cdot \nabla V_d
\,-\,i\,\frac{4 |V_d|^2 \psi_2}{(1+|V_d|^2)^2},
\\
\\
{\mathbb N}_2 [\psi] &= & \frac{1}{V_d} \nabla \psi \cdot \nabla V_d \, O(\psi_2)
\,+\,  O\Big(   \big| |V_d|^2-1\big| +|\psi_2|\Big) \big|\nabla \psi \cdot \nabla \psi \big|
\,+\,  i\, O\big( | e^{-\psi_2} -1+\psi_2|\big),
\\
\\
{\mathbb M}_2 [\psi] &=&
2\,\ve\,\frac{\partial V_d}{\partial x_2}\frac{ 2|V_d|^2 (1-e^{-2\psi_2})}{\, (1+|V_d|^2 e^{-2 \psi_2})(1+|V_d|^2)\,V_d\,}
\,+\,2\,i\,\ve\,\frac{\partial \psi}{\partial x_2}\frac{ 1-2|V_d|^2e^{-2\psi_2}}{\, 1+|V_d|^2 e^{-2 \psi_2}\,\,},
\\
\\
{E} &=& -\frac{{\mathbb E}_1}{\,i V_d \,}.
\end{eqnarray*}
Recall that $\psi=\psi_1+i \psi_2$. Then setting $ z=\vec{e}_1 + y$, we have
for $z \in \R^{2}_{+}$
\begin{align}
\label{L0vortexpairs}
{\mathbb L}_2 [\psi]= \begin{pmatrix}
  \Delta \psi_1 + O(e^{-|y|} ) |\nabla \psi|
\\
\\
\Delta \psi_2 -\frac{4 |V_d|^2}{(1+|V_d|^2)^2} \psi_2 +O(e^{-|y|}) \nabla \psi_2
\end{pmatrix},\qquad\qquad\qquad
\\
\nonumber
\\
\label{N1nvortexpairs}
{\mathbb M}_2 [\psi]=\begin{pmatrix}
  O\big(e^{-|y|}  |\nabla \psi \cdot \nabla \psi |+ |\psi_2|^2 \frac{1}{(1+|y|)^2}
  +|\psi_2| \frac{1}{1+|y|} |\nabla \psi|\big)
\\
\\
  O\big(e^{-|y|}  |\nabla \psi \cdot \nabla \psi |+ |\psi_2| |\nabla \psi \cdot \nabla \psi |
   +|\psi_2|^2\big)
\end{pmatrix},\qquad
\\
\nonumber
\\
\label{N2vortexpairs}
{\mathbb N}_2 [\psi]=\begin{pmatrix}
  O(\ve |\frac{\partial \psi_2}{\partial y_2}|
\\
\\
  O(\ve |\frac{\partial \psi_1}{\partial y_2}|
\end{pmatrix}.\qquad\qquad\qquad\qquad\qquad\qquad
\end{align}
Let us remark that the explicit form of all the linear and nonlinear terms will be very useful for later analysis.
We can also write problem (\ref{orignalproblem-Q2T2-vortexpairs}) as the above in a similar way,
but the nonlinear operator ${\mathbb M}_2$ has the form
\begin{align*}
{\mathbb M}_2 [\psi]
\,=\,&
2\, \ve\,\frac{\partial V_d}{\partial x_2}\frac{ 2|V_d|^2 (1-e^{-2\psi_2})}{\, (1+|V_d|^2 e^{-2 \psi_2})(1+|V_d|^2)\,V_d\,}
\,+\,2\,i\, \ve\,\frac{\partial \psi}{\partial x_2}\frac{ 1-2|V_d|^2e^{-2\psi_2}}{\, 1+|V_d|^2 e^{-2 \psi_2}\,\,},
\\
&\,-\,2\, \kappa\ve\,\frac{\partial V_d}{\partial x_2}\frac{ 2|V_d|^2 (1-e^{-2\psi_2})}{\, (1+|V_d|^2 e^{-2 \psi_2})(1+|V_d|^2)\,V_d\,}
\,-\,2\,i\, \kappa\ve\,\frac{\partial \psi}{\partial x_2}\frac{ 1-2|V_d|^2e^{-2\psi_2}}{\, 1+|V_d|^2 e^{-2 \psi_2}\,\,}.
\end{align*}

\medskip
Let us fix two small positive numbers $ 0< \gamma <1,\,  0<\varrho <1$.
Recall that $ \phi= iV_d \psi,\, \psi= \psi_1+ i \psi_2$.
Denote $ \ell_j= |z-\vec{e}_j|$, and define
\begin{eqnarray}
\| \psi \|_{*} &=& \sum_{j=1}^2 \| \phi \|_{C^{2, \gamma} ( \ell_j <2)} + \sum_{j=1}^2 \| \phi_j\|_{C^{1, \gamma} (\ell_j <3)}  \label{5.1n} \\
& &  + \sum_{j=1}^2 \Bigg[ \| \ell_j^\varrho \psi_1 \|_{L^\infty (\ell_j >2)} +\| \ell_j^{1+\varrho} \nabla \psi_1 \|_{L^\infty (\ell_j >2)} \Bigg]
 \nonumber \\
& &  + \sum_{j=1}^2 \Bigg[ \| \ell_j^{1+\varrho} \psi_2 \|_{L^\infty (\ell_j >2)} +\| \ell_j^{2+\varrho} \nabla \psi_2 \|_{L^\infty (\ell_j >2)} \Bigg]
 \nonumber
\end{eqnarray}
\begin{equation}
\| h\|_{**}= \sum_{j=1}^2 \| iV_d h \|_{C^{0, \gamma} (\ell_j<3)} + \sum_{j=1}^2
\Big[\,\| \ell_j^{2+\varrho} h_1\|_{L^\infty (\ell_j >2)} + \| \ell_j^{1+\varrho} h_2 \|_{L^\infty (\ell_j >2)}\,\Big].
\end{equation}

We remark that the choices of these norms are motivated by the expressions of (\ref{N1nvortexpairs})-(\ref{N2vortexpairs}).
A direct application of (\ref{S0(Vd)}), (\ref{S1(Vd)}), and (\ref{S2(Vd)}) yields
the decay estimates for the errors
$$
\tilde{E}_1\,=\, -\frac{{\mathbb E}_1}{\,i V_d \,},
\qquad
\tilde{E}_2\,=\, -\frac{{\mathbb E}_2}{\,i V_d \,},
$$
which will be stated in the lemma.
\begin{lemma}
\label{l5.1vortexpairs}
It holds that for $  z \in \big(B_{2R} (\vec{e}_1) \cup B_{2R} (\vec{e}_2)\big)^c$
\begin{align}
\label{Ere0}
\Big|\, {\mathrm {Re}} (\tilde{E}_1),\,  {\mathrm {Re}} (\tilde{E}_2)\,\Big|
\,\leq\,
\frac{C \ve^{1-\varrho}}{ (1+ |z-\vec{e}_1|)^{3}}
\,+\,
\frac{C \ve^{1-\varrho}}{ (1+ |z-\vec{e}_2|)^{3}},
\\
\nonumber
\\
\label{Eim0}
\Big|\, {\mathrm {Im}}(\tilde{E}_1),\, {\mathrm {Re}} (\tilde{E}_2)\,\Big|
\,\leq\, \frac{C \ve^{1-\varrho}}{ (1+ |z-\vec{e}_1|)^{1+\varrho}}
\,+\,
\frac{C \ve^{1-\varrho}}{ (1+ |z-\vec{e}_2|)^{1+\varrho}},
\end{align}
where $\varrho \in (0,1)$ is a constant. Moreover,
\begin{align}
\big\|\, \tilde{E}_1,\, \tilde{E}_2\,\big\|_{**}
\,\leq\,
C\ve^{1-\varrho}.
\end{align}
\end{lemma}

\subsection{Projected Nonlinear Problem}


%

Let $\tilde{\eta}$ be defined as in (\ref{cutoff}) and $R>0$ be a fixed large positive number.
By defining
\begin{equation}
\label{Zd}
Z_d:= \frac{\partial V_d}{\partial d}
\Bigg[\,\tilde{\eta} \Big(\frac{|z-\vec{e}_1|}{R}\Big)+ \tilde{\eta}\Big(\frac{|z-\vec{e}_2|}{R}\Big)\,\Bigg],
\end{equation}
we consider the full nonlinear {\em projected problem}
\begin{align}
\begin{aligned}
\label{nonlinearprojection-pairs}
{\mathcal L}[\psi]
\,+\,
{\mathcal N}[\psi]
\,+\,
{\mathcal M}[\psi]
\,=\,
{\mathcal E}
\,+\,
c\,Z_d,
\qquad\qquad
\\
\mbox{Re} \Big(\int_{ \R^2 } \bar{\phi}\, Z_d\Big)\,=\,0,
\qquad
\psi\mbox{ satisfies the symmetry (\ref{psisymvirtexpairs})},
\end{aligned}
\end{align}
where we have denoted that
\begin{align*}
{\mathcal L}\,=\,{\mathbb L}_1,
\quad
{\mathcal N}\,=\,{\mathbb N}_1,
\quad
{\mathcal M}\,=\,{\mathbb M}_1
\quad\mbox{in}\quad {\mathfrak S}_1,
\\
{\mathcal L}\,=\,{\mathbb L}_2,
\quad
{\mathcal N}\,=\,{\mathbb N}_2,
\quad
{\mathcal M}\,=\,{\mathbb M}_2
\quad\mbox{in}\quad {\mathfrak S}_2.
\end{align*}
Note that in the above we have used the relation $\phi\,=\,iV_d\psi$ in ${\mathfrak S}_1$.
It is easy to show that $Z_d$ possesses the symmetries in (\ref{symmetryvorxpairs}).
Here is the resolution theory.

\begin{proposition}
\label{Proposition4.1}
There exists a constant $C$, depending on $\gamma, \varrho$
only such that for all $\ve$ sufficiently small, $d$ large,
the following holds: there exists a  unique solution $\psi_\ve$
to (\ref{nonlinearprojection-pairs}) and $ \psi_\ve$ satisfies
\begin{equation}
\| \psi_\ve \|_{*} \leq C \ve^{1-\varrho}.
\end{equation}
Furthermore, $ \psi_{\ve, d}$ is continuous in $d$.
\qed
\end{proposition}
The proof of this proposition is similar of that for Proposition \ref{Proposition5.1}.

\subsection{Reduced Problem}\label{reductionprocedurepairs}

From Proposition  \ref{Proposition4.1},
we deduce the existence of  a solution $(\phi, c)$ to (\ref{nonlinearprojection-pairs}).
To find a real solution to (\ref{orignalproblem-Q2-vortexpairs}), we shall choose suitable $d$
such that $c$ is zero.
This can be realized by the standard reduction procedure in this subsection.

\medskip
On $\R_{+}^2, $ we recall that $z= \vec{e}_1 +y$ and
$$
V_d= w (y) e^{-i \theta_{\vec{e}_2}} \big(1+ O(e^{-d/2-|y|/2})\big),
$$
then
\begin{align*}
Z_d \,=\,& \frac{\partial V_d}{\partial d}
\,=\, - \frac{\partial w}{\partial y_1} e^{-i \theta_{\vec{e}_2}} \big(1+ O(e^{-d/2-|y|/2})\big)
\,-\, i w \frac{\partial \theta_{\vec{e}_2}}{\partial d}   e^{-i \theta_{\vec{e}_2}} \big(1+ O(e^{-d/2-|y|/2})\big),
\end{align*}
where ${\partial \theta_{\vec{e}_2}}/{\partial d}$ is small in the sense
\begin{equation*}
\frac{\partial \theta_{\vec{e}_2}}{\partial d}
\,=\,\frac{\partial \theta_{\vec{e}_2}}{\partial x_1}
\,=\, O(\frac{1}{d})
\,=\, O(\ve).
\end{equation*}
In the polar coordinates, we have $w= \rho (r) e^{i \theta}$
\[
\frac{\partial w}{\partial y_1}= \big( \rho' \cos \theta - i \frac{\rho}{r} \sin \theta\big) e^{i \theta},
\qquad
\frac{\partial w}{\partial y_2}= \big( \rho' \sin \theta + i \frac{\rho}{r} \cos \theta\big) e^{i \theta}.
\]

\medskip
Multiplying (\ref{nonlinearprojection-pairs}) by $ \frac{1}{(1+|V_d|^2)^2}\overline{Z_d} $ and integrating, we obtain
\begin{align*}
c\, \mbox{Re} \Bigg(\,\int_{\R^2} \frac{1}{(1+|V_d|^2)^2}\,Z_d\, {\overline {Z_d}}\,\Bigg)
\,=\,
&\,-\,\mbox{Re} \Bigg(\,\int_{\R^2} \frac{1}{(1+|V_d|^2)^2}\,  \overline{Z_d}\, {\mathcal E}\,\Bigg)
\\
&\,+\,\mbox{Re} \Bigg(\int_{\R^2} \frac{1}{(1+|V_d|^2)^2}\, \overline{Z_d}\,
\Big({\mathcal L}[\psi]+{\mathcal N}[\psi]+{\mathcal M}[\psi]\Big)\Bigg)
\\
\,=\,
&\,-\,\mbox{Re} \Bigg(\,\int_{\R^2} \frac{1}{(1+|V_d|^2)^2}\,  \overline{Z_d}\, {\mathbb E}_d\,\Bigg)
\\
&\,+\,\mbox{Re} \Bigg(\int_{\R^2} \frac{1}{(1+|V_d|^2)^2}\, \overline{Z_d}\,
\Big( {\mathbb L}_1[\phi]+{\mathbb N}_1[\phi]+{\mathbb M}_1[\phi]\Big)\Bigg).
\end{align*}
Using Proposition \ref{Proposition4.1} and the expression in (\ref{N1nvortexpairs}), we deduce that
\begin{equation}
 \mbox{Re} \Bigg(\int_{\R^2} \overline{Z_d}\, {\mathbb N}_1 [\phi]\Bigg)
 \,=\, o(\ve).
\end{equation}
On the other hand, integration by parts, we have
\begin{align*}
\mbox{Re} \Bigg(\int_{\R^2} \frac{1}{(1+|V_d|^2)^2} \, \overline{Z_d}\, {\mathbb L}_1 [\phi]\Bigg)
 \,=\,
 \mbox{Re} \Bigg(\int_{\R^2} \frac{1}{(1+|V_d|^2)^2} \, \overline{\phi}\, {\mathbb L}_1[Z_d]\Bigg).
\end{align*}
Let us observe that
\begin{equation}
\frac{\partial}{\partial d} {\mathbb S}_0 [V_d]
\,=\, {\mathbb L}_1 \Big[\frac{\partial V_d}{\partial d}\Big]
\,=\, {\mathbb L}_1 [Z_d] = o(\ve),
\end{equation}
and thus by Proposition \ref{Proposition4.1}
\begin{equation}
\mbox{Re} \Bigg(\int_{\R^2} \frac{1}{(1+|V_d|^2)^2}
\,\overline{\phi} \,{\mathbb L}_1 [Z_d]\Bigg)
\,=\, o(\ve).
\end{equation}

It remains to estimate the following integral
\begin{align*}
-2\,  \mbox{Re} \Bigg(\int_{\R^2_{+} } \frac{1}{(1+|V_d|^2)^2}
\,\overline{Z_d}\,{\mathbb E}_d [V_d] \Bigg)
\,=\,&\,
\mbox{Re} \Bigg(\int_{\R^2} \frac{1}{(1+|V_d|^2)^2} \, \overline{Z_d}\, {\mathbb S}_1 [V_d] \Bigg)
\\
\,=\,&\,2\,  \mbox{Re} \Bigg(\int_{\R^2_{+} } \frac{1}{(1+|V_d|^2)^2}\, \overline{Z_d}
\, {\mathbb S}_0 [V_d]\Bigg)
\\
&\,+\, 2\,\ve\,\mbox{Re} \Bigg( i \int_{\R^2_{+} } \frac{1}{(1+|V_d|^2)^2}\,  \overline{Z_d}
\,  {\mathbb Q}_2 [V_d]\Bigg),
\end{align*}
where $ {\mathbb S}_0[V_d]$ and $ {\mathbb Q}_2[V_d]$ are defined in (\ref{S0(Vd)}) and (\ref{Q2(Vd)}).

\medskip
The expression in (\ref{Q2(Vd)}) gives that
\begin{align}\label{IntegralQ2Pairs}
\begin{aligned}
&\ve\,\mbox{Re} \Bigg(\int_{\R^2} \frac{1}{(1+|V_d|^2)^2} \,{\mathbb Q}_2 [V_d]\, \overline{Z_d}\Bigg)
\\
\,=\,&
\,-\,\ve\, \mbox{Re}\Bigg( i \int_{\R^2} \frac{1-|w|^2}{(1+|w|^2)^3} \frac{\partial w}{\partial y_2} \overline{\frac{\partial w}{\partial y_1}}\Bigg)
\,+\,o(\ve)
\\
\,=\,&
\,-\,\ve\, \mbox{Re}\Bigg( i \int_{\R^2} \frac{1-|\rho|^2}{(1+|\rho|^2)^3}
\big( \rho' \sin \theta + i \frac{\rho}{r} \cos \theta\big)
\big( \rho' \cos \theta + i \frac{\rho}{r} \sin \theta\big)
\Bigg)
\,+\,o(\ve)
\\
\,=\,&\,\ve\, \int_{\R^2} \frac{(1-\rho^2)\rho \rho'}{r (1+\rho^2)^3 }
\,=\, -\frac{\pi\ve}{4} \,+\,o(\ve).
\end{aligned}
\end{align}
On the other hand, using the estimate (\ref{S0(Vd)}), we have
\begin{eqnarray*}
 & &\mbox{Re} \Bigg(  \int_{\R_{+}^2} \frac{1}{(1+|V_d|^2)^2}  \,{\mathbb S}_0 [V_d]\, \overline{Z_d}\Bigg)
 \\
&=& -\,\mbox{Re} \Bigg(\int_{\R_{+}^2} \frac{1}{(1+|V_d|^2)^2} \,{\mathbb S}_0 [ V_d] \,
\overline{\frac{\partial w}{\partial y_1}}\, e^{i \theta_{\vec{e}_2}} \Bigg)
\,+\, o(\ve)
\\
&= & -\,\mbox{Re} \Bigg( \int_{\R_{+}^2} \frac{1}{(1+|w|^2)^2}
\Bigg[ \frac{2 (|w|^2 -1)}{|w|^2+1} i \nabla w \cdot \nabla \theta_{\vec{e}_2}
- \frac{1-|w|^2}{1+|w|^2} \frac{w}{|y+\vec{e}_1-\vec{e}_2|^2}\Bigg]
\\
&&\qquad\qquad\times e^{-i \theta} \big( \rho' \cos \theta + i \frac{\rho}{r} \sin \theta\big)\Bigg)
\,+\,o(\ve).
\end{eqnarray*}
Note that in $\R^2_{+}$
\begin{equation}
\nabla \theta_{\vec{e}_2} \approx \frac{\vec{e}_2}{2d},
\quad
\nabla w \cdot \nabla  \theta_{\vec{e}_2} \approx \frac{e^{i \theta}}{2d}
\big[ \rho' \sin \theta +i \frac{\rho}{r} \cos \theta\big].
\end{equation}
Hence
\begin{eqnarray*}
& & \mbox{Re} \Bigg(\int_{\R^2} \frac{1}{(1+|V_d|^2)^2} \,{\mathbb S}_0 [ V_d]\, \overline{Z_d}\Bigg)
\\
& = & -\,\mbox{Re} \Bigg(\int_{\R^2}  \frac{ (|w|^2-1)}{(|w|^2+1)^3} \frac{i}{d}
\Big[\,\rho' \sin\theta + i \frac{\rho}{r} \cos \theta\,\Big]
\Big[\, \rho' \cos \theta +i \frac{\rho}{r} \sin \theta\,\Big]\Bigg)
\,+\,o(\ve)
\\
&= & \frac{1}{d} \int_{\R^2} \frac{ (\rho^2-1)}{(\rho^2+1)^3} \frac{\rho' \rho}{r}
\,+\,o(\ve)
\\
& = & \frac{\pi}{d} \int_0^\infty \frac{\rho^2-1}{(\rho^2+1)^3} (\rho^2)'
\,+\,o(\ve)
\,=\, -\,\frac{\pi}{4d}\,+\,o(\ve).
\end{eqnarray*}

\medskip
Combining all estimates together, we obtain the following equation
\begin{equation}\label{reducealgebravorpairQ2}
c (d)= c_0 \Big[\,-\frac{\pi}{4 d}  +\frac{\pi\ve}{4}\,\Big]  + o(\ve),
\end{equation}
where $ o(\ve)$ is continuous function of $d$
(which is a consequence of continuity of $\phi$ in $d$) and $c_0 \not =0$.
By simple mean-value theorem, we can find a zero of $c(d)$.
This completes the reduction procedure
for problem (\ref{orignalproblem-Q2-vortexpairs}),
which also finishes the proof of the existence of vortex pairs of (\ref{Wavemaps})
in Theorem \ref{theorem-vortexpairs}.

\medskip
For the reduction procedure of problem (\ref{orignalproblem-Q2T2-vortexpairs}),
we shall estimate
\begin{align*}
c_\ve\, \mbox{Re} \Bigg(\,\int_{\R^2} \frac{1}{(1+|V_d|^2)^2}\,Z_d\, Z_d\,\Bigg)
\,=\,&\,\mbox{Re} \Bigg(\,\int_{\R^2} \frac{1}{(1+|V_d|^2)^2}\,  \overline{Z_d}
\, {\mathbb S}_2 [V_d]\,\Bigg)
\\
&\,+\,
\mbox{Re} \Bigg(\,\int_{\R^2} \frac{1}{(1+|V_d|^2)^2}\, \overline{Z_d}
\, {\mathbb L}_1 [\phi]\,\Bigg)
\\
&\,+\,\mbox{Re} \Bigg(\int_{\R^2} \frac{1}{(1+|V_d|^2)^2}\, \overline{Z_d}
\,\big( {\mathbb N}_1 [\phi]+{\mathbb M}_1 [\phi]\big)\Bigg),
\end{align*}
where ${\mathbb S}_2\,=\,{\mathbb S}_0\,+\,\ve\,{\mathbb Q}_2\,+\,\kappa\,\ve\,{\mathbb T}_2$
is defined in (\ref{S0QT}).
The expression in (\ref{T2(Vd)}) gives that
\begin{align}\label{IntegralT2Pairs}
\begin{aligned}
&\kappa\,\ve\,\mbox{Re} \Bigg( i \int_{\R^2} \frac{1}{(1+|V_d|^2)^2} \,{\mathbb T}_2 [V_d]\, \overline{Z_d}\Bigg)
\\
\,=\,&\,\kappa\,\ve\, \mbox{Re}\Bigg( i \int_{\R^2} \frac{1}{(1+|w|^2)^2} \frac{\partial w}{\partial y_2} \overline{\frac{\partial w}{\partial y_1}}\Bigg)
\,+\,o(\ve)
\\
\,=\,&\,\kappa\,\ve\, \mbox{Re}\Bigg( i \int_{\R^2} \frac{1}{(1+|\rho|^2)^2}
\big( \rho' \sin \theta + i \frac{\rho}{r} \cos \theta\big)
\big( \rho' \cos \theta + i \frac{\rho}{r} \sin \theta\big)
\Bigg)
\,+\,o(\ve)
\\
\,=\,&\, -\,\kappa\,\ve\, \int_{\R^2} \frac{\rho \rho'}{r (1+\rho^2)^2 }
\,=\, -\,\kappa\,\ve\,\frac{\pi}{2} \,+\,o(\ve).
\end{aligned}
\end{align}
Whence after similar computation in the above, we need consider the equation
\begin{equation}\label{reducealgebravorpairQ2T2}
c(d)= c_0 \Big[\,-\frac{\pi}{4 d}+\frac{\pi\ve}{4}-\frac{\ve\,\kappa\,\pi}{2}\,\Big]  + o(\ve),
\end{equation}
where $ o(\ve)$ is continuous function of $d$
(which is a consequence of continuity of $\phi$ in $d$) and $c_0 \not =0$.
If $1-2\kappa\neq 0$,
by simple mean-value theorem, we can find a zero of $c(d)$.
This completes the reduction procedure
for problem (\ref{orignalproblem-Q2T2-vortexpairs}),
which also finishes the proof of the existence of vortex pairs of (\ref{Ishimoriproblem})
in Theorem \ref{theorem-vortexpairs}.

\begin{remark}\label{remarkLocationandveloPairs}
Recall the parameters $\ve$, $\kappa$ given in (\ref{vePairs}), (\ref{kvePairs}).
By the relation in (\ref{reducealgebravorpairQ2}), we have that, for problem (\ref{orignalproblem-Q2-vortexpairs})
(i.e. (\ref{Wavemaps})), the vortex and antivortex pair
undergoes Kelvin motion
when its speed is sufficiently small so that the vortices are widely separated with the relation
$$
\frac{1}{d}\sim \frac{2c}{\sqrt{1-c^2}}.
$$
It is also that from (\ref{reducealgebravorpairQ2T2}),
for problem (\ref{orignalproblem-Q2T2-vortexpairs}) (i.e. (\ref{Ishimoriproblem})), there holds
$$
\frac{1}{d}\sim \frac{2c}{\sqrt{1-c^2}}-\frac{2\omega}{\sqrt{1-c^2}}.
$$
\end{remark}

\section{Vortex helices}\label{section5}
\setcounter{equation}{0}

In this section, we will construct solutions with vortices to
problem (\ref{orignalproblem-Q2T2-vortexhelices}),
which will provide vortex rings to problem (\ref{Ishimoriproblem}).
We omit the similar arguments for the existence of solution to problem (\ref{orignalproblem-Q2-vortexhelices})
(i.e. the existence of vortex rings to (\ref{Wavemaps})).

\subsection{Approximate solutions}\label{subsection51}

\subsubsection{First Approximate Solution and its Error}

Recall the vortex solutions $w^{+}$ and $w^{-}$ defined in (\ref{8}).
For each fixed $d :=\frac{\hat{d}}{\ve} $ with $\hat{d} \in [  {1}/{ 100},  \    100  ] $,
we define the {\bf first approximate solution}
\begin{equation}
v_0(z):= w^{+} (z-\vec{e}_1) w^{-} (z-\vec{e}_2),
\end{equation}
where ${\vec e}_1\,=\,(d,0)$ and ${\vec e}_2\,=\,(-d,0)$.
A simple computation shows that
\begin{equation}
|v_0(z)| \to 1 \ \mbox{as} \ |z| \to +\infty.
\end{equation}
In fact,  for $|z| >>d$, we have
\begin{eqnarray*}
 v_0&\approx & e^{i \theta_{\vec{e}_1} - i \theta_{\vec{e}_2}}
 \\
& \approx & \frac{ (x_1^2 -d^2 + x_2^2 +2i x_2 d)}{\, \sqrt{ (x_1-d)^2 +x_2^2}\, \sqrt{ (x_1+d)^2 +x_2^2} }
\\
& \approx & 1 + O\Big(\frac{d^2}{|z|^2} +\frac{d}{|z|} \Big).
\end{eqnarray*}
Here $\theta_{\xi}$ denotes the angle argument around $ \xi$ and $ \ell_\xi=|z-\xi|$.   It is easy to see that
\begin{equation}\label{prederivates}
\nabla \ell_{\xi}= \frac{1}{|z-\xi|} (z-\xi),
\quad
\nabla \theta_{\xi}= \frac{1}{|z-\xi|^2} (z-\xi)^\perp,
\quad
|\nabla \theta_\xi |=\frac{1}{|z- \xi|},
\end{equation}
where  we denote $ z^\perp= (-x_2, x_1)$.

\medskip
By our construction and the properties of $\rho$, we have
\begin{equation}
v_0= e^{i\theta_{\vec{e}_1} -i\theta_{\vec{e}_2}}
\Big[\,1+ O\big(e^{- \min  (|z-\vec{e}_1|, |z-\vec{e}_2|)}\big)\,\Big].
\end{equation}
We divide the region ${\mathfrak S}$ into two parts:
\begin{align}
\begin{aligned}
{\mathfrak S}_+=\Big\{\,z: x_1>0\,\Big\},
\qquad
{\mathfrak S}_-=\Big\{\,z: x_1<0\,\Big\}.
\end{aligned}
\end{align}
By our symmetry assumptions, we just need to consider the region ${\mathfrak S}_+$.
In the region ${\mathfrak S}_+$, we consider two cases.
Firstly, if $ |z-\vec{e}_2| > \frac{d}{2}+ \frac{|z-\vec{e}_1|}{2}$, then we have
\begin{eqnarray*}
 v_0&=& w(z-\vec{e}_1) e^{-i \theta_{\vec{e}_2}} \Big[\,1+ O(e^{-d/2 -|z-d\vec{e}_1|/2 })\,\Big],
 \\
 |v_0| &=&  \rho (|z-\vec{e}_1|) \Big[\,1+ O(e^{-d/2-|z-\vec{e}_1|/2})\,\Big].
\end{eqnarray*}
Secondly, if $ |z-\vec{e}_2| < \frac{d}{2}+ \frac{|z-\vec{e}_1|}{2}$,
in the rest of this paper, we often use the  translated variable as follows
\begin{align}\label{translatedcoordinates}
z\,=\,\vec{e}_1 + y.
\end{align}

\medskip
The term ${\mathbb H}_1[v_0]$ obeys the following asymptotic behavior
\begin{align*}
\Omega_1&\,:=\,{\mathbb H}_1[v_0]\,=\,\frac{1}{x_1}\frac{\partial v_0}{\partial x_1}
\\
&\,=\, \frac{1}{d+y_1}\frac{\partial}{\partial y_1} \Bigg[ w (y) e^{-i \theta_{\vec{e}_2}} + O( e^{-d/2-|y|/2})\Bigg]
\\
&\,=\,\Big(\frac{1}{d}-\frac{y_1}{d(d+y_1)}\Big) \Bigg[ \frac{\partial w}{\partial y_1}
-i w \frac{\partial \theta_{\vec{e}_2}}{\partial y_1} \Bigg]
e^{-i \theta_{\vec{e}_2}}
\,+\, O(\ve e^{-d/2-|y|/2} ),
\end{align*}
where we have
\begin{equation}
\frac{\partial }{\partial y_1} \theta_{\vec{e}_2}
\,=\, \frac{ -x_2}{(y_1+2d)^2+y_2^2}
\,=\, O\big(d^{-2}\big)
\,=\, O(\ve^2).
\end{equation}
Whence, by (\ref{prederivates}),
the term $\frac{1}{d} \frac{\partial w}{\partial y_1}$ in $\Omega_1$ has a singularity in the form
\begin{align}\label{singularity1}
\frac{i}{d}w \frac{\partial \theta_{\vec{e}_1}}{\partial y_1}/(\,i v_0\,)
\,\sim\,-\frac{1}{d}\frac{y_2}{|y|^2}.
\end{align}

\medskip
In summary, by recalling (\ref{S0(Vd)}), (\ref{Q2(Vd)}), and (\ref{T2(Vd)}),
we have obtained for $ z \in {\mathfrak S}_+$ with $z= \vec{e}_1 +y$
\begin{align}
\begin{aligned}\label{S4(v0)}
{\mathbb S}_4 [v_0]&\,=\,
{\mathbb S}_0[v_0]
\,+\,\ve\,|\log\ve|\, \,{\mathbb Q}_2[v_0]
\,+\,\kappa\,\ve\,|\log\ve|\, \,{\mathbb T}_2[v_0]
\,+\,{\mathbb H}_1[v_0]
\\
&\,=\,  e^{-i \theta_{\vec{e}_2}}
\Bigg[
\frac{2 (\rho^2-1)}{\rho^2+1} i \nabla w \cdot \nabla \theta_{\vec{e}_2}
+ \frac{\rho^2 -1}{\rho^2+1 } \frac{w}{|y+\vec{e}_1-\vec{e}_2|^2}
\Bigg]
\\
&\qquad\,+\,
i \ve|\log\ve|e^{-i \theta_{\vec{e}_2}}
\Bigg[\frac{\partial w}{\partial y_2}
-iw \frac{\partial \theta_{\vec{e}_2}}{\partial y_2}  + O(e^{-d/2-|y|/2} )
\Bigg]\frac{1-\rho^2}{1+\rho^2}
\\
&\qquad\,-\,
i\kappa \ve|\log\ve|e^{-i \theta_{\vec{e}_2}}
\Bigg[  \frac{\partial w}{\partial y_2}
-iw \frac{\partial \theta_{\vec{e}_2}}{\partial y_2}  + O(e^{-d/2-|y|/2} )
\Bigg]
\,+\,\Omega_1.
\end{aligned}
\end{align}
The similar (and almost identical) estimates  also hold in the region $z \in {\mathfrak S}_-$.
A direct application of (\ref{S4(v0)})
yields  the following decay estimates for
the error
\begin{align}
\begin{aligned}
{\tilde E}_4 \,=\, \frac{\,
{\mathbb S}_0 [v_0]
\,+\,\ve\,|\log\ve|\, \,{\mathbb Q}_2 [v_0]
\,+\,\kappa\,\ve\,|\log\ve|\, \,{\mathbb T}_2 [v_0]
\,}{i v_0}.
\end{aligned}
\end{align}

\begin{lemma}\label{l5.1vortexhelix}
It holds that for $  z \in \big(B_2 (\vec{e}_1) \cup B_2 (\vec{e}_2)\big)^c$
\begin{equation}
\label{Ere}
\Big|{\mathrm {Re}} (\tilde{E}_4)\Big|
\,\leq\, \frac{C \ve^{1-\varrho}}{ (1+ |z-\vec{e}_1|)^{3}}
\,+\, \frac{C \ve^{1-\varrho}}{ (1+ |z-\vec{e}_2|)^{3}},
\end{equation}
\begin{equation}
\label{Eim}
\Big|{\mathrm {Im}}(\tilde{E}_4)\Big|
\,\leq\, \frac{C \ve^{1-\varrho}}{ (1+ |z-\vec{e}_1|)^{1+\varrho}}
\,+\, \frac{C \ve^{1-\varrho}}{ (1+ |z-\vec{e}_2|)^{1+\varrho}},
\end{equation}
where $\varrho \in (0,1)$ is a constant.
\end{lemma}

\medskip

\noindent
{\bf Proof:} The estimates for the term ${\mathbb S}_0 [v_0]\frac{ 1}{\,i v_0\,}$
follows from (\ref{S4(v0)}).
We just need to  estimate the terms
$$
\ve\,|\log\ve|\, \,{\mathbb Q}_2[v_0]\frac{ 1}{\,i v_0\,},
\qquad
\kappa\,\ve\,|\log\ve|\, \,{\mathbb T}_2[v_0]\frac{ 1}{\,i v_0\,},
\qquad
{\mathbb H}_3[v_0],
$$
whose definition are given in (\ref{S0QT}) and (\ref{H-operators}).

\medskip
Let us compute for $z \in {\mathfrak S}_+$
\begin{eqnarray*}
 \ve i \frac{\partial v_0}{\partial x_2}\,\Big/i v_0
&=&
\ve \frac{ \partial \big(\rho (z-\vec{e}_1) e^{i \theta_{ \vec{e}_1}- \theta_{\vec{e}_2}}\big)}
{\partial x_2}\frac{1}{  v_0}
\,+\,O\big(\ve e^{- |z-\vec{e}_1|} +\ve e^{-|z-\vec{e}_2|}\big)
\\
&=& \ve \frac{\partial \rho (|z-\vec{e}_1|)}{\partial x_2}\frac{1}{\rho (|z- \vec{e}_1|)}
\, +\, i \ve  \frac{ \partial (\theta_{\vec{e}_1} - \theta_{\vec{e}_2})}{\partial x_2}
\,+\,O\big( \ve e^{-|z- \vec{e}_1|}+ \ve e^{-|z-\vec{e}_2|}\big)
\\
& =& O\big(\ve e^{- |z-\vec{e}_1|} +\ve e^{-|z-\vec{e}_2|}\big)
\,+\, i \ve \Bigg[  \frac{ x_1-d}{(x_1-d)^2+x_2^2} -\frac{x_1+d}{(x_1+d)^2+x_2^2} \Bigg].
\end{eqnarray*}
The estimate (\ref{Ere}) then follows.
Let us notice that for $ z \in {\mathfrak S}_+, \  |z-\vec{e}_1| <d$,
\begin{equation*}
\ve \Big| \frac{ x_1-d}{(x_1-d)^2+x_2^2} -\frac{x_1+d}{(x_1+d)^2+x_2^2} \Big|
\, \leq\, C \ve \frac{1}{ (1+|z-\vec{e}_1|)}
\,\leq\, C \ve^{1-\varrho} \frac{1}{ (1+|z-\vec{e}_1|)^{1+\varrho}}.
\end{equation*}
On the other hand, for $ z \in {\mathfrak S}_+,\  |z- \vec{e}_1| >d$, we then have
\begin{equation*}
\ve \Big| \frac{ x_1-d}{(x_1-d)^2+x_2^2} -\frac{x_1+d}{(x_1+d)^2+x_2^2} \Big|
\,\leq\, C \frac{\ve }{ (1+|z-\vec{e}_1|)^2}
\,\leq\, C  \frac{\ve^{1-\varrho}}{ (1+|z-\vec{e}_1|)^{1+\varrho}}.
\end{equation*}
Thus (\ref{Eim}) is proved.
\qed

\subsubsection{Further Improvement of the Approximation}

As we promised in previous subsection, we now define a new correction $\varphi_d$ in the phase term.
The phase function $ \varphi_d$ will be decomposed into two parts:
singular part and regular part.
Then we will define an improved approximation and estimate
its error by substituting to problem (\ref{orignalproblem-Q2T2-vortexhelices}).

\medskip
To cancel the singularities in (\ref{singularity1}),
we want to find a function $\Phi(y_1,y_2)$ by solving the problem in the translated coordinates $(y_1,y_2)$
\begin{align}\label{green}
\frac{\partial^2\Phi}{\partial y_1^2}
\,+\,\frac{\partial^2\Phi}{\partial y_2^2}
\,=\,
\frac{y_2}{d\gamma^2|y|^2}
\quad\mbox{in }\R^2.
\end{align}
In fact, we can solve this problem by separation of variables and then obtain
\begin{align}
\begin{aligned}\label{correctionpahse}
\Phi(y_1,y_2)&\,=\,\frac{1}{\,4d\gamma^{2}\,}y_2\log|y|^2.
\end{aligned}
\end{align}

\medskip
Let $\chi$ be a smooth cut-off function in a way such that
$ \chi (z)= 1$ for $ z  \in B_{\frac{d}{10}} ( \vec{e}_1)$
and $ \chi=0$  for $ z \in \big(B_{\frac{d}{5}} (\vec{e}_1)\big)^c $.
Let $ \varphi_{d} (z)= \varphi_s (z) + \varphi_{r} (z)$,
where the singular part is defined by
\begin{equation}
\begin{aligned}\label{varphis}
\varphi_{s} (z)&\,:=\,
\chi (z)\frac{1}{\,4 d\,\gamma^2} x_2 \log \frac{ |z- \vec{e}_1|^2}{ \,|z-\vec{e}_2|^2\,}.
\end{aligned}
\end{equation}
While by recalling the operators ${\mathbb H}_1$ in (\ref{H-operators}),
we find the regular part $ \varphi_{r} (z) $ by solving the problem
\begin{align*}
\Big[\Delta + {\mathbb H}_1 \Big]\varphi_r
\,=\,&
  - \Big[\Delta + {\mathbb H}_1 \Big]
\big( \theta_{ \vec{e}_1} -\theta_{\vec{e}_2} + \varphi_s\big)\quad{\mbox{in}}\quad\R.
\end{align*}
Note that the function $ \varphi_{s}$ is continuous but $\nabla \varphi_s$ is not.
The singularity of $\varphi_s$ comes from its derivatives.

\medskip
By simple computations, we see that for $ z \in B_{\frac{d}{10}} (\vec{e}_1)$,
\begin{align*}
&  \Big[\Delta + {\mathbb H}_1\Big]
\big( \theta_{ \vec{e}_1} -\theta_{\vec{e}_2} + \varphi_s\big)
\\
&\quad\,=\, \frac{ 4 x_2 (x_1-d)}{ \,\big((x_1-d)^2+ x_2^2\big)\,\big( (x_1+d)^2+ x_2^2\big)\,}
\,+\,
\frac{1}{x_1} \frac{ x_2 (x_1^2-x_2^2-d^2)}{\, \big((x_1-d)^2+ x_2^2\big)\,\big( (x_1+d)^2+ x_2^2\big)\,}
\\
&\quad\,=\, O(d^{-2})= O(\ve^2).
\end{align*}
For   $ z \in \big(B_{\frac{d}{10}} (\vec{e}_1)\big)^c$, it is easy to see that we also get $O(\ve^2)$.
In fact for $z \in \big(B_{\frac{d}{5}} (\vec{e}_1)\big)^c$, $\varphi_s=0$ and
\begin{align*}
\Big[\Delta + {\mathbb H}_1\Big]
\big( \theta_{ \vec{e}_1} -\theta_{\vec{e}_2} + \varphi_s \big)
\,=\,\frac{ -4 x_2 }{\, \big((x_1-d)^2+ x_2^2\big)\big( (x_1+d)^2+ x_2^2\big)\,}.
\end{align*}
Going back to the original variable $(r, {{\check s}_3 })$ in (\ref{coordinatesrs})
and letting $ \hat{\varphi} (r,{{\check s}_3 })= \varphi_r (z)$ we see that
\begin{equation}
\Bigg| \Delta_{r, {{\check s}_3 }} \hat{\varphi}
\,+\, {\mathbb H}_1[\hat{\varphi}]
\Bigg|
\,\leq\,
\frac{C}{\, \big( \, \sqrt{1 + r^2+ {{\check s}_3 }^2} \, \big)^3\,}.
\end{equation}
Thus we can choose $ \varphi_r $ such
that $ \hat{\varphi}= O\big(\frac{1}{ \sqrt{1+ r^2+ |{{\check s}_3 }|^2}} \big) $.
The regular term $\varphi_r$ is $ C^{1}$ in the original variable $(r, {{\check s}_3 })$.

\medskip
We observe also that by our definition, the function
\begin{equation}
\label{tildevar}
\tilde{\varphi}:= \theta_{ \vec{e}_1} - \theta_{\vec{e}_2} + \varphi_s + \varphi_r,
\end{equation}
satisfies
\begin{equation}
\label{fact4}
\Big[\Delta + {\mathbb H}_1\Big]{\tilde{\varphi}}\, =\,0
\quad\mbox{on}\quad\R.
\end{equation}
From the decomposition of $ \varphi_d$,
we see that the singular term contains $ x_2 \log  |z-\vec{e}_1|$
which becomes dominant when we calculate the speed.

\medskip
Finally,  we define an improved approximation
\begin{equation}\label{Vd-helix}
V_d (z):= w^{+} (z-\vec{e}_1) w^{-} (z-\vec{e}_2) e^{ i \varphi_d }.
\end{equation}

\subsection{Error Estimates}\label{subsection4.3}

Let
\[
\tilde{\rho} = \rho ( |z-\vec{e}_1|) \rho ( |z-\vec{e}_2|),
\quad
\tilde{\varphi} = \theta_{\vec{e}_1} -\theta_{\vec{e}_2} + \varphi_d.
\]
So $ V_d= \tilde{\rho} e^{i \tilde{\varphi} }$.  Note that for $ x_1>0$, we have
\begin{equation}
\rho (|z-\vec{e}_2|)= 1+ O(e^{- \frac{d}{2} - \frac{1}{2} |z-\vec{e}_1|}).
\end{equation}
Since the error between $1$ and $ \rho (|z-\vec{e}_2|)$ is exponentially small, we may ignore $\rho (|z-\vec{e}_2|)$ in the computations below.
We shall check that $V_d$ is a good approximate solution in the sense that
it satisfy the conditions in (\ref{symmetryandboundary})
and has a small error.

\medskip
Let us start to compute the errors:
\begin{equation}
\Delta V_d= \Biggl[  \Delta \tilde{\rho} - |\nabla \tilde{\varphi} |^2 \tilde{\rho }
+ 2 i\nabla \tilde{\rho} \cdot \nabla \tilde{\varphi}
-i \tilde{\rho}\,{\mathbb S}_2[{\tilde\varphi}]  \Biggl] e^{i \tilde{\varphi} }.
\end{equation}
Here we have used the fact (\ref{fact4}).
We continue to compute other terms:
\[
\nabla V_d \cdot \nabla V_d\,=\,
\Biggl[ |\nabla \tilde{\rho} |^2- \tilde{\rho}^2 |\nabla \tilde{\varphi}|^2+ 2i \tilde{\rho} \nabla \tilde{\rho} \cdot \nabla \tilde{\varphi} \Biggl]
e^{2 i \tilde{\varphi}}.
\]
We then obtain that
\begin{align}
\begin{aligned} \label{errornew}
{\mathbb S}_0 [V_d]
&\,=\, e^{i \tilde{\varphi}}
\Biggl[ -i \tilde{\rho}\,{\mathbb S}_0[{\tilde\varphi}]
\,+\, \frac{ \tilde{\rho} (\tilde{\rho}^2-1)}{\tilde{\rho}^2+1}\,
\big( |\nabla \tilde{\varphi}|^2 - |\nabla \theta|^2 \big)
\\
&\qquad\qquad\,+\,2 i \frac{1-\tilde{\rho}^2}{1+\tilde{\rho}^2} \nabla \tilde{\rho} \cdot \nabla \tilde{\varphi}
\,+\,O\Big(e^{-\frac{d}{2}-\frac{|z-\vec{e}_1|}{2}}\Big)
\Biggl].
\end{aligned}
\end{align}

In a small neighborhood of ${\vec{e}_1}$, we write
\begin{align*}
\begin{aligned}
{\mathbb H}_1[V_d]
&\,=\,
\frac{1}{x_1}\frac{\partial V_d}{\partial x_1}
\\
&\,=\,
\frac{1}{x_1}\frac{\partial{\tilde\rho}}{x_1}e^{i{\tilde\varphi}}
\,+\,i{\tilde\rho}\,{\mathbb H}_1{\tilde\varphi}
\\
&\,=\,i{\tilde\rho}\,{\mathbb H}_1{\tilde\varphi}
\,+\,O(\ve).
\end{aligned}
\end{align*}
The estimates for the terms
\begin{align*}
\kappa\, \ve |\log \ve |{\mathbb T}_2[V_d]
\,=\,&\,
i\kappa\, \ve |\log \ve | \frac{\partial V_d}{\partial x_2},
\\
 \ve |\log \ve |{\mathbb Q}_2[V_d]
\,=\,&\,
i\, \ve |\log \ve |
\frac{1-|V_d|^2}{1+|V_d|^2}
\frac{\partial V_d}{\partial x_2}
\end{align*}
are same as before.

\medskip
The total error is
\begin{align}
\begin{aligned}\label{S4(Vd)}
{\mathbb E}_d&\,=\, \frac{ \tilde{\rho} (\tilde{\rho}^2-1)}{\tilde{\rho}^2+1}\,
\big( |\nabla \tilde{\varphi}|^2 - |\nabla \theta|^2 \big)
\,+\,2 i \frac{1-\tilde{\rho}^2}{1+\tilde{\rho}^2} \nabla \tilde{\rho} \cdot \nabla \tilde{\varphi}
\,+\,O\Big(e^{-\frac{d}{2}-\frac{|y|}{2}}\Big)
\\
&\qquad\,+\,i \ve|\log\ve|  \Bigg[ \frac{\partial w}{\partial y_2} -i w \frac{\partial \theta_{\vec{e}_2}}{\partial y_2} \Bigg]
\frac{1-\rho^2}{1+\rho^2}\,e^{-i \theta_{\vec{e}_2}}
\\
&\qquad\,+\,i \ve|\log\ve|  \Bigg[ \frac{\partial w}{\partial y_2} -i w \frac{\partial \theta_{\vec{e}_2}}{\partial y_2} \Bigg]
\,e^{-i \theta_{\vec{e}_2}}
\,+\, O(\ve|\log\ve| e^{-d/2-|y|/2} ).
\end{aligned}
\end{align}
Here we have used the relation in (\ref{fact4}).

\medskip
Setting $ z\,=\,\vec{e}_1+ y$, we then have
\begin{equation}
\nabla \varphi_s= - \frac{1}{\,2 d\gamma^2\,}  \log d \, \nabla y_2 \,+\, O(\ve \log |y| ),
\qquad
\nabla \varphi_r =O(\ve).
\end{equation}
Thus
\begin{eqnarray*}
 |\nabla \tilde{\varphi}|^2- |\nabla \theta|^2
 &=& O(|\nabla \varphi_s \cdot \nabla \theta |)+ |\nabla \varphi_s|^2 + O(\ve |\nabla \theta |),
\\
 \nabla \tilde{\rho} \cdot \nabla \tilde{\varphi}
 &= & O( \ve |\log \ve | \rho'   )+ O(\ve  \rho' \log |y| ).
\end{eqnarray*}
These asymptotic expression will play an important role in the reduction part.

\subsection{Setting up of the Problem}\label{settingup}

Now we introduce the set-up of the reduction procedure.
We look for solutions to  (\ref{orignalproblem-Q2-vortexhelices}) and problem (\ref{orignalproblem-Q2T2-vortexhelices})
with boundary condition in (\ref{boundary}) in the form
\begin{equation}
\label{perturbationdecompositionhelix}
 u (z)= \eta ( V_d + i V_d \psi) \,+\, (1-\eta)\, V_d\, e^{i \psi},
\end{equation}
where $\eta$ is a function such that
\begin{equation}
\label{eta}
 \eta= \tilde{\eta} (|z-\vec{e}_1|)+ \tilde{\eta} (|z-\vec{e}_2|),
\end{equation}
and $ \tilde{\eta} (s)= 1$ for $ s \leq 1$ and $ \tilde{\eta} (s)=0$ for $s \geq 2$.

\medskip
The conditions imposed on $u$ in (\ref{boundary}) and (\ref{symmetryvortexhelix})
can be transmitted to the symmetry on $\psi$
\begin{align}\label{psisym}
\begin{aligned}
\psi(z)= \psi (- \bar{z} ),
\qquad
\psi (z)=-\overline{\psi ({\bar z})},
\qquad
\frac{\partial\psi}{\partial x_1}(0,x_2)\,=\,0.
\end{aligned}
\end{align}
More precisely,
for $\psi=\psi_1+i\psi_2$, there hold the conditions
\begin{align}
\begin{aligned}\label{psisymmetryandboundary}
\psi_1(x_1,x_2)=\psi_1(-x_1,x_2),
&\qquad
\psi_1(x_1,x_2)=-\psi_1(x_1,-x_2),
\\
\psi_2(x_1,x_2)=\psi_2(-x_1,x_2),
&\qquad
\psi_2(x_1,x_2)=\psi_2(x_1,-x_2),
\\
\frac{\partial \psi_1}{\partial x_1}(0,x_2)=0,
&\qquad
\frac{\partial \psi_2}{\partial x_1}(0,x_2)=0.
\end{aligned}
\end{align}
This symmetry will be important in solving the linear problems
in that it excludes all but one kernel.
We may write $\psi=\psi_1+ i \psi_2$ with $\psi_1, \psi_2$ real-valued and then set
\begin{equation}
u\,=\, V_d +\phi,
\quad
\phi\,=\, \eta i V_d \psi + (1-\eta) V_d (e^{i \psi} -1).
\end{equation}
In the sequel, we will derive the explicit local form of the equation for the perturbation term $\psi$.

\medskip
Let $R>1$ be a fixed constant.
In the inner region
$$
{\mathfrak S}_1\,=\,\big\{\,z\,|\, z \in B_{9R} (\vec{e}_1) \cup B_{9R} (\vec{e}_2)\,\big\},
$$
we have
\begin{equation}
u\,=\,V_d \,+\,\phi,
\end{equation}
and the equation for $\phi$ becomes
\begin{equation}
{\mathbb L}_1 [\phi] \,+\, {\mathbb N}_1 [\phi]\,=\, {\mathbb E}_d.
\end{equation}
In the above, we have denoted linear operator by
\begin{align*}
\begin{aligned}
{\mathbb L}_d [\phi]\,=\,&\,
\Delta \phi
\,+\,\frac{1}{x_1}\frac{\partial\phi}{\partial x_1}
\,-\,\frac{4 {\bar V}_d}{1+|V_d|^2} \nabla V_d \cdot \nabla \phi
\,-\,\frac{2 \bar{\phi}}{1+|V_d|^2} \nabla V_d \cdot \nabla V_d
\\
&\,+\, \frac{\,2 {\bar V}_d (V_d \bar{\phi} + {\bar V}_d \phi)\,}{(1+|V_d|^2)^2} \nabla V_d \cdot \nabla V_d
\,+\, F'(V_d) \phi.
\end{aligned}
\end{align*}
The nonlinear operator and the error are
\begin{align}
\begin{aligned}
\label{Nphi-helices}
{\mathbb N}_1[\phi]\,=\,&\, F(V_d+\phi)
\,-\, F(V_d)
\,-\,F'(V_d) \phi
\,+\, O\Big((1+|\phi|) |\nabla \phi|^2\Big)
\\
&\,+\,i\kappa\,\delta\, |\log\ve|\ve \frac{\partial \phi}{\partial x_2}
\,+\,i\kappa\,\delta\, |\log\ve|\ve \frac{\partial \phi}{\partial x_2},
\end{aligned}
\end{align}
and
\begin{equation}
{\mathbb E}_d\,=\,-{\mathbb S}_4 [V_d].
\end{equation}
In the above, we have use the definition of $F$ in (\ref{F}).

\medskip
In the outer region
$$
{\mathfrak S}_2\,=\,\Big\{\, z\,|\, z \in \big(B_{4R} (\vec{e}_1) \cup B_{4R} (\vec{e}_2)\big)^c\,\Big\},
$$
we have $ u= V_d e^{i \psi}$.
By simple computations we obtain
\begin{align*}
\frac{ {\mathbb S} [V_d e^{i \psi}]}{iV_de^{i\psi}} \,=\, &\,
\Delta \psi
\,+\,\frac{1}{x_1}\frac{\partial\psi}{\partial x_1}
\,+\, 2 \frac{ 1-|V_d|^2 + |V_d|^2 (e^{-2\psi_2} -1)}{ V_d (1+|V_d|^2 e^{-2\psi_2})} \nabla V_d \cdot \nabla \psi
\\
&\,+\, \frac{1}{iV_d} \frac{2 |V_d|^2 {\bar{V}}_d (e^{-2\psi_2} -1)}{ (1+|V_d|^2)(1+|V_d|^2 e^{-2\psi_2})} \nabla V_d \cdot \nabla V_d
\,-\,i \Big(\frac{2 |V_d|^2}{1+|V_d|^2 e^{-2\psi_2}} -1\Big) \nabla \psi \cdot \nabla \psi,
\\
&\,-\,i  \frac{ 2|V_d|^2 (1-e^{-2\psi_2})}{ (1+|V_d|^2 e^{-2 \psi_2})(1+|V_d|^2)}
\,+\,i\delta |\log\ve| \ve  \frac{\partial \psi}{\partial x_2}
\\
&\,+\,i\kappa\delta |\log\ve| \ve\frac{\,1-|V_d|^2e^{-2\psi_2}\,}{1+|V_d|^2e^{-2\psi_2}}
\frac{\partial \psi}{\partial x_2}
\,+\, \frac{{\mathbb E}_d}{iV_d} \, .
\end{align*}
We can also write the problem as an equation of $\psi=\psi_1+i\psi_2$
\begin{equation}
{\mathbb L}_2 [\psi]
\,+\,{\mathbb N}_2 [\psi]
\,+\,{\mathbb M}_2 [\psi]
\,=\,E,
\end{equation}
with conditions in (\ref{psisymmetryandboundary}).
In the above, we have denoted
\begin{align*}
\begin{aligned}
 {\mathbb L}_1 [\psi] \,=\,& \Delta \psi
\,+\,\frac{1}{x_1}\frac{\partial\psi}{\partial x_1}
\,+\, \frac{2 (1-|V_d|^2)}{1+|V_d|^2} \nabla V_d \cdot \nabla \psi
\,+\, \frac{ 4i {\bar V}_d^2 \psi_2}{(1+|V_d|^2)^2} \nabla V_d \cdot \nabla V_d
\\
&\,-\,i \frac{4 |V_d|^2 \psi_2}{(1+|V_d|^2)^2},
\end{aligned}
\end{align*}
\begin{align*}
 {\mathbb N}_2 [\psi] \,=\,&\, \frac{1}{V_d} \nabla \psi \cdot \nabla V_d\,O(\psi_2)
\,+\,  O\big(   \big| |V_d|^2-1\big| +|\psi_2|\big) |\nabla \psi \cdot \nabla \psi |
\, +\,  i O( | e^{-\psi_2} -1+\psi_2|)
\\
&\,+\,i\,\kappa\,\delta\, |\log\ve|\ve \frac{\partial \psi}{\partial x_2}\,.
\end{align*}
There error term has the form
\begin{align*}
E \,=\, -\frac{{\mathbb E}_d}{i V_d }.
\end{align*}
Recall that $\psi=\psi_1+i \psi_2$. Then setting $ z=\vec{e}_1 + y$, we have
for $z \in \R^{2}_{+}$
\begin{align}
\label{L0}
{\mathbb L}_2 [\psi]
\,=\,
\left[
\begin{array}{cc}
\Delta \psi_1
\,+\,\frac{1}{x_1}\frac{\partial\psi_1}{\partial x_1}
\,+\, O(e^{-|y|} ) |\nabla \psi|
\\
\\
\Delta \psi_2
\,+\,\frac{1}{x_1}\frac{\partial\psi_2}{\partial x_1}
\, -\,\frac{4 |V_d|^2}{(1+|V_d|^2)^2} \psi_2
\,+\,O\big(e^{-|y|}\big) \nabla \psi_2
\end{array}
\right],
\end{align}

\begin{align}
\label{N1n}
{\mathbb N}_2 [\psi]
\,=\,
\left[
\begin{array}{cc}
  O\Big(e^{-|y|}  |\nabla \psi \cdot \nabla \psi |
\,+\, |\psi_2|^2 \frac{1}{(1+|y|)^2}
\,+\,|\psi_2| \frac{1}{1+|y|} |\nabla \psi|\Big)
\\
\\
  O\Big(e^{-|y|}  |\nabla \psi \cdot \nabla \psi |
\,+\, |\psi_2| |\nabla \psi \cdot \nabla \psi |
\,+\,|\psi_2|^2\Big)
\end{array}
\right],
\end{align}

\begin{align}
\label{N2}
{\mathbb N}_2 [\psi]
\,=\,
\left[
\begin{array}{cc}
  O(\ve |\frac{\partial \psi_2}{\partial y_2}|)
\\
\\
  O(\ve |\frac{\partial \psi_1}{\partial y_2}|)
\end{array}
\right].
\end{align}
Let us remark that the explicit form of all the linear and nonlinear terms
will be very useful for later analysis in resolution theory.

\medskip
Let
\begin{equation}
{\mathbb E}_4\,=\, {\mathbb S}_4 [ V_d],
\quad
{{\tilde{\mathbb E}}_4}= \frac{ {\mathbb E}_4}{i V_d}.
\end{equation}
Based on the form of the errors, we need to use suitable norms.
Let us fix two  positive numbers
$$
p >13,\quad  0<\varrho <1.
$$
Recall that $ \phi= iV_d \psi,\,\, \psi= \psi_1+ i \psi_2$.
Denote $ \ell_1= |z- \vec{e}_1|$ and $ \ell_2= |z- \vec{e}_2|$, and define
\begin{align*}
\| h\|_{**}\,=\,  \| iV_d h \|_{L^p ({\mathfrak S}_1)}
\,+\,  \sum_{j=1}^2\Big[\,\| \ell_j^{2+\varrho} h_1\|_{L^\infty ({\mathfrak S}_2)}
+ \| \ell_j^{1+\varrho} h_2 \|_{L^\infty ({\mathfrak S}_2)} \,\Big],
\end{align*}
\begin{align*}
\begin{aligned}
\| \psi \|_{*}&\,=\,
\| \phi \|_{ W^{2,p} ( {\mathfrak S}_1)}
\,+\,\sum_{j=1}^2\Bigg[ \| \ell_j^\varrho \psi_1 \|_{L^\infty ({\mathfrak S}_2)} +\| \ell_j^{1+\varrho} \nabla \psi_1 \|_{L^\infty ({\mathfrak S}_2)} \Bigg]
\\
&\qquad
\,+\, \Bigg[ \| \ell_j^{1+\varrho} \psi_2 \|_{L^\infty ({\mathfrak S}_2)} +\| \ell_j^{2+\varrho} \nabla \psi_2 \|_{L^\infty ({\mathfrak S}_2)} \Bigg].
\end{aligned}
\end{align*}
In the above,
\begin{align}
\begin{aligned}
{\mathfrak S}_1\,=\,\{\,z\in{\mathfrak{S}}: |z-{\vec e}_1|<3\,\, \mbox{or}\,\, |z-{\vec e}_2|<3\,\},
\\
{\mathfrak S}_2\,=\,\{\,z\in{\mathfrak{S}}: |z-{\vec e}_1|>2\,\, \mbox{and}\,\, |z-{\vec e}_2|>2\,\}.
\end{aligned}
\end{align}
We remark that we use the norm $ L^p_{loc} $ (or $W^{2,p}_{loc}$) in the inner part
due to the fact that  the error term  contains terms
like $ \ve \log |z-\vec{e}_1|$ which is not $L^\infty$-bounded.

\medskip
Using the norms defined above, we can have the following error estimates.
\begin{lemma}\label{l8.1m}
It holds that for $  z \in {\mathfrak S}_2$
\begin{equation}
\label{Ere-8}
\big| {\mathrm{Re}}(\tilde{{\mathbb E}}_4)\big|
\,\leq\,
\frac{C \ve^{1-\varrho}}{ (1+ |z-\vec{e}_1|)^{3 }}
\,+\,
\frac{C \ve^{1-\varrho}}{ (1+ |z-\vec{e}_2|)^{3 }},
\end{equation}
\begin{equation}
\label{Eim-8}
\big| {\mathrm {Im}}(\tilde{{\mathbb E} }_4)\big|
\,\leq\,
\frac{C \ve^{1-\varrho}}{ (1+ |z-\vec{e}_1|)^{1+\varrho}}
\,+\,
\frac{C \ve^{1-\varrho}}{ (1+ |z-\vec{e}_2|)^{1+\varrho}},
\end{equation}
and also
\begin{equation}
\| {\mathbb E}_d \|_{L^p ({\mathfrak S}_1)} \,\leq\, C \ve |\log \ve |,
\end{equation}
where $\varrho \in (0,1)$ is a constant. As a consequence, there holds
\begin{equation}
\|\tilde{{\mathbb E}}_4\|_{**} \,\leq\, C \ve^{1-\varrho}.
\end{equation}
\qed
\end{lemma}

\subsection{Projected Linear and Nonlinear Problems}\label{projectedproblem}

Let
\begin{eqnarray*}
L_0 [\psi] &=& \frac{\partial^2 \psi}{\partial x_1^2}
\,+\, \frac{\partial^2\psi}{\partial x_2^2}
\,+\,\frac{1}{x_1}\frac{\partial\psi}{\partial x_1}
\\
& \ &
\,+\,  \frac{2 (1-|V_d|^2)}{1+|V_d|^2} \nabla V_d \cdot \nabla \psi
\,+\, \frac{ 4i {\bar V}_d^2 \psi_2}{(1+|V_d|^2)^2} \nabla V_d \cdot \nabla V_d
\,-\,i \frac{4 |V_d|^2 \psi_2}{(1+|V_d|^2)^2},
\end{eqnarray*}
and the co-kernel
\begin{equation}
\label{Zd1-8}
Z_d:= \frac{\partial V_d}{\partial d} \Bigg[\tilde{\eta} \Big(\frac{|z-\vec{e}_1|}{R}\Big)
+ \tilde{\eta} \Big(\frac{|z-\vec{e}_2|}{R}\Big)\Bigg],
\end{equation}
where $\tilde{\eta}$ is defined at (\ref{eta}) before.
Then $Z_d$ satisfies the symmetry (\ref{symmetryvortexhelix}).

\medskip
As the first step of finite dimensional reduction, we  need to consider the following linear problem
\begin{align}
\label{linear2-8}
\begin{aligned}
L_0 (\psi)\,=\, h  \ \ \mbox{in} \ {\mathfrak S},
\qquad
\psi \ \mbox{satisfies the symmetry (\ref{psisym})},
\\
\\
 {\mathrm{Re}} \Big(\int_{ {\mathfrak S}} \bar{\phi} Z_d\,\mathrm{d}x\Big)\,=\,0
\quad \mbox{for }\phi\,=\, iV_d \psi \mbox{ in } B_1 (\vec{e}_1)\cup B_1 (\vec{e}_2).
\end{aligned}
\end{align}
We have the following {\em a priori estimates}.
\begin{lemma}
\label{8.1n}
There exists a constant $C$, depending on $\varrho$ only such that for all $\ve$ sufficiently small, $d \sim \frac{1}{\ve}$, and any solution of (\ref{linear2-8}), it holds
\begin{equation}
\| \psi \|_{*} \leq C \| h \|_{**}.
\end{equation}
\end{lemma}

\noindent
{\bf Proof:} The proof is similar as in Lemma 5.1 in \cite{linwei}.
 Suppose that there exists a sequence of $\ve= \ve_n \to 0$,
functions $\psi^n, h_n$ which satisfy (\ref{linear2-8}) with
\begin{equation*}
\| \psi^n \|_{*}=1,\quad \| h_n \|_{**} =o(\ve).
\end{equation*}
We will derive a contradiction by careful analysis of the estimates.

\medskip
We derive {\em inner estimates} first.
We have the symmetries and boundary conditions for $\psi_1$
and $\psi_2$ in (\ref{psisymmetryandboundary}).
Whence we may just need to consider the region
\[
\Sigma_+\,=\,\big\{\, x_1>0,\,\, 0<x_2<\lambda\pi/(\ve\gamma)\,\big\}.
\]
Then we have
\begin{equation*}
 \mbox{Re} \Bigg(\,\int_{ \R^2 } \bar{\phi}_n Z_d\,\Bigg)
 \,=\,
 2  \mbox{Re} \Bigg(\,\int_{ {\Sigma_+}} \bar{\phi}_n Z_d\,\Bigg)
 \,=\,0.
\end{equation*}
Let $ z \in \Sigma_+,\,\, z= \vec{e}_1 +y$ and $ \tilde{\phi}_n (y)= \phi_n (z)$. Then as $ n \to +\infty$,
\begin{equation*}
 V_d\,=\, w^{+} (y) e^{-i \theta_{\vec{e}_2}} \big(1+ O(e^{-d/2})\big)
 \,=\, -w^{+} (y) + o(1).
\end{equation*}
Since $ \|\psi^n \|_{*}=1$, we may take a limit so that $ \tilde{\phi}_n \to \phi_0$ in $ \R^2_{loc}$, where $ \phi_0$ satisfies
\begin{equation*}
{\mathbb L}_0 [\phi_0]=0
\end{equation*}
where ${\mathbb L}_0$ is defined by (\ref{L1n}). Observe that $ \phi_0$
satisfies the decay estimate (\ref{psidecay}) because of our assumption on $\psi^n$.
By Theorem \ref{non1}, we have
\[ \phi_0= c_1 \frac{\partial w}{\partial y_1}+ c_2 \frac{\partial w}{\partial y_2}.\]
Observe that $\phi_0$ inherits the symmetries of $\phi$ and hence $ \phi_0 = \overline{\phi_0 (\bar{z})}$.
(The other symmetry is not preserved under the transformation $z= \vec{e}_1 +y$.)
But certainly  $\frac{\partial w}{\partial y_2}$  does not enjoy the above symmetry.
Hence $ \phi_0=  c_1\frac{\partial w}{\partial y_1}$.
On the other hand, taking a limit of the orthogonality condition
$$
\mbox{Re} \Bigg(\int_{ {\Sigma_+}} \bar{\phi}_n Z_d\Bigg)\,=\,0,
$$
we obtain
$$
\mathrm{Re}\Bigg(\,\int_{\R^2} \bar{\phi_0} \frac{\partial w}{\partial y_1}\,\Bigg)\,=\,0.
$$
This implies that $ c_1=0$ and hence we have
\begin{equation*}
\phi_n \to 0 \ \ \mbox{in} \ \R^2_{loc}
\end{equation*}
which implies that for any fixed $R>0$,
\begin{equation}\label{inner}
\sum_{j=1}^2 \Big(\,\| \phi_1\|_{L^p ( \ell_j <R)}
\,+\,\| \phi_2\|_{L^p ( \ell_j <R)}
\,+\, \| \nabla \phi_1\|_{L^p ( \ell_j <R)}
\,+\,\| \nabla \phi_2\|_{L^p ( \ell_j <R)}\Big)
\,=\,o(1).
\end{equation}
We use  the $L^p$-estimates in the inner part $ \{ | z-\vec{e}_1| <R \}$.
By choosing $p>N $ large we obtain the embedding $ W^{2, p}_{loc}$ into $ C^{1, \alpha}_{loc}$
for any $\alpha \in (0,1)$.

\medskip
Next we shall derive {\em outer estimates}.
Let $\tilde{\eta}$ be a cut-off function such that $ \tilde{\eta} (s)=1$ for $ s\leq 1$
and $ \tilde{\eta} (s)=0$ for $ s>2$.
We consider the new function
\[ \tilde{\psi} \,=\, \psi\,\chi (z),  \ \ \mbox{where} \
\chi (z)\,=\,
1\,-\,
\tilde{\eta} \Big(\frac{|z-\vec{e}_1|}{4} \Big)
\,-\,\tilde{\eta} \Big(\frac{|z-\vec{e}_2|}{4}\Big).
\]

Using the explicit forms of $L_0$ in (\ref{L0}), the first equation becomes
\begin{align}
\begin{aligned}\label{5.30}
&\frac{\partial^2 {\tilde\psi}_1}{\partial x_1^2}
\,+\,\frac{1}{x_1}\frac{\partial{\tilde\psi}_1}{\partial x_1}
\,+\,\frac{\partial^2{\tilde\psi}_1}{\partial x_2^2}
\,=\, O(e^{-|y|}) |\nabla \psi|
\,+\, O( \nabla \chi \nabla \psi_1)
\,+\, O(\psi_1 \Delta \chi)
\,+\, h_1 \chi.
\end{aligned}
\end{align}
On the region $\Sigma\backslash (B_4 (\vec{e}_1) \cup  B_{4} (\vec{e}_2))$,
we have the conditions in (\ref{psisymmetryandboundary}).
Moreover we have
\[
|\nabla \chi \cdot \nabla \psi|
\,=\, o(1)\big(|z-\vec{e}_1|^2 +|z-\vec{e}_2|^2\big)^{-\frac{\sigma+2}{2}}.
\]

For the outer part estimates,  we use the following new barrier function
\begin{equation*}
B(z):\,=\, B_1(z)\,+\, B_2(z),
\end{equation*}
where we have denoted
\begin{align}
\begin{aligned}\nonumber
B_1(z)\,=\,\ell_1^{\beta}\,x_2^{\nu}+\ell_2^{\beta}\,x_2^{\nu},
\qquad
B_2 (z)\,=\,C_1\big(1+|z|^2\big)^{-\varrho/2}.
\end{aligned}
\end{align}
In the above $\beta+\nu=-\varrho,\,\, 0<\varrho <\nu<1$.

\medskip
Now, we do the computations for $B_1$
\begin{equation*}
\frac{\partial^2 B_1}{\partial x_1^2}
\,+\,\frac{\partial^2 B_1}{\partial x_2^2}
\,\leq\,
-C\,\Big(\,{\ell}_1^{-2-\varrho}\,+\,{\ell}_2^{-2-\varrho}\,\Big).
\end{equation*}
\begin{equation*}
\frac{\partial^2 B_2}{\partial x_1^2}
\,+\,\frac{\partial^2 B_2}{\partial x_2^2}
\frac{1}{x_1} \frac{\partial B_2}{\partial x_1}
\,\leq\,
-C\,C_1\,\big(\,1+|z|^2\,\big)^{-1-\varrho/2}.
 \end{equation*}
where $C_1$ depends only on $\beta$ and $\nu$.
On the other hand, there holds
\begin{equation*}
\frac{1}{x_1} \frac{\partial B_1}{\partial x_1}
\,\leq\,
C \frac{x_2^\nu }{x_1}
 \Big[\, {\ell}_1^{\beta -2} (x_1-d)+ {\ell}_2^{\beta-2} (x_1+d)\,\Big],
\end{equation*}
Thus  for $\ell_1 < c_\varrho d$ or $ \ell_2 < c_\varrho d$,
where $c_\varrho $ is small, we have
\begin{equation*}
\frac{1}{x_1} \frac{\partial B_1}{\partial x_1}
\,\leq\,
C_2\, c_\varrho\big(\,{\ell}_1^{-2-\varrho}\,+\,{\ell}_2^{-2-\varrho}\,\big).
\end{equation*}
where $C_2$ depends only on $\beta$ and $\nu$.
For ${\ell_1} > c_\varrho d$, there holds
\begin{equation*}
\frac{1}{x_1} \frac{\partial B_1}{\partial x_1} \leq C  \big(1+|z|^2\big)^{-1-{\varrho}/{2}}.
\end{equation*}
By choosing $C_1$ large, we finally have
\begin{equation*}
\frac{\partial^2 B}{\partial x_1^2}
\,+\,\frac{\partial^2 B}{\partial x_2^2}
\,+\,\frac{1}{x_1}\frac{\partial B}{\partial x_1}
\,\leq\,
\,-\,C \big(\,\ell_1^{-2-\varrho} +\ell_2^{-2-\varrho}\,\big).
\end{equation*}

By comparison principle on the set $\Sigma\backslash ( B_4 (\vec{e}_1) \cup  B_{4} (-\vec{e}_1))$, we get that
\begin{equation}
\label{out1}
 |\tilde{\psi}_1 | \leq C B  ( \| h\|_{**}+ o(1)),
 \quad
 \forall  \ z \in \Sigma \backslash ( B_4 (\vec{e}_1) \cup  B_{4} (-\vec{e}_1)).
\end{equation}
Elliptic estimates then give
\begin{equation}
\label{out2}
\sum_{j=1}^2  \big\| \ell_j^{1+\sigma}   |\nabla \tilde{\psi}_1 |\big\|_{L^\infty (\ell_j >4)}
\,\leq\,
C ( \| h\|_{**}+ o(1)).
\end{equation}

To estimate $\psi_2$, we perform the same cut-off and now the second equation becomes
\begin{align}
\begin{aligned}\label{5.31}
&\frac{\partial^2 {\tilde\psi}_2}{\partial x_1^2}
\,+\,\frac{1}{x_1}\frac{\partial{\tilde\psi}_2}{\partial x_1}
\,+\,\frac{\partial^2{\tilde\psi}_2}{\partial x_2^2}
 -\frac{4 |V_d|^2}{(1+|V_d|^2)^2} \tilde{\psi}_2
\\
&\,=\,O(\frac{1}{1+|y|}) \nabla \psi_1
\,+\, O(e^{-|y|}) \nabla \psi_2
\,+\, O( \nabla \chi \nabla \psi_2)
\,+\, O(\Delta \psi_2)
\,+\,h_2 \chi \, .
\end{aligned}
\end{align}
Note that we also have the conditions in (\ref{psisymmetryandboundary}). Since for $ z \in \Sigma \backslash ( B_4 (\vec{e}_1) \cup  B_{4} (-\vec{e}_1))$, there holds
$$ \frac{4|V_d|^2}{(1+|V_d|^2)} \geq \frac{1}{4},$$
by standard elliptic estimates we have
\begin{align}
\begin{aligned}\label{out4}
\|\psi_2 \|_{L^\infty (\ell_j >4)}
\,\leq\,
C (\| \psi_2 \|_{L^\infty (\ell_j=4)}) (1+\| \psi\|_{*}) \| h \|_{**}
(1+|z-\vec{e}_1| + |z-d\vec{e}_1|)^{-1-\varrho},
\\
\\
|\nabla \psi_2|
\,\leq\,
C (\| \psi_2 \|_{L^\infty (\ell_j=R)}) (1+\| \psi\|_{*}) \| h \|_{**}
(1+|z-\vec{e}_1| + |z-d\vec{e}_1|)^{-2-\varrho}.
\end{aligned}
\end{align}

\medskip
Combining both inner and outer estimates in (\ref{inner}), (\ref{out1})-(\ref{out2})
and (\ref{out4}),
we obtain that $ \|\psi\|_{*}=o(1)$, which is a contradiction.
\qed

\medskip
We consider the full nonlinear {\em projected problem}
\begin{align}
\begin{aligned}
\label{nonlinearprojection-helix}
{\mathcal L}[\psi]
\,+\,
{\mathcal N}[\psi]
\,+\,
{\mathcal M}[\psi]
\,=\,
{\mathcal E}
\,+\,
c\,Z_d,
\qquad\qquad
\\
\mbox{Re} \Big(\int_{ \R^2 } \bar{\phi}\, Z_d\Big)\,=\,0,
\qquad
\psi\mbox{ satisfies the symmetry (\ref{psisymmetryandboundary})},
\end{aligned}
\end{align}
where we have denoted that
\begin{align*}
{\mathcal L}\,=\,{\mathbb L}_1,
\quad
{\mathcal N}\,=\,{\mathbb N}_1,
\quad
{\mathcal M}\,=\,{\mathbb M}_1
\quad\mbox{in}\quad {\mathfrak S}_1,
\\
{\mathcal L}\,=\,{\mathbb L}_2,
\quad
{\mathcal N}\,=\,{\mathbb N}_2,
\quad
{\mathcal M}\,=\,{\mathbb M}_2
\quad\mbox{in}\quad {\mathfrak S}_2.
\end{align*}
Note that in the above we have used the relation $\phi\,=\,iV_d\psi$ in ${\mathfrak S}_1$.
By contraction mapping theorem, we make a conclusion by the following the resolution theory
\begin{proposition}
\label{Proposition5.1}
There exists a constant $C$, depending on $p, \varrho$
only such that for all $\ve$ sufficiently small, $d$ large,
the following holds: there exists a  unique solution $\psi_{\ve,d}$
to (\ref{nonlinearprojection-helix}) and $ \psi_{\ve,d}$ satisfies
\begin{equation}
\| \psi_\ve \|_{*} \leq C \ve^{1-\varrho}.
\end{equation}
Furthermore, $ \psi_{\ve, d}$ is continuous in $d$.
\qed
\end{proposition}

\medskip
Here, we do not give the proof to the last proposition.
The reader can refer to the arguments in \cite{linwei} for details.

\subsection{Reduced Problem}\label{reductionprocedurehelix}

From Proposition  \ref{Proposition5.1},
we deduce the existence of  a solution $(\psi, c)$ to (\ref{nonlinearprojection-helix}).
To find a real solution to problem (\ref{orignalproblem-Q2T2-vortexhelices}),
we shall choose suitable $d$ such that $c$ is zero.
This can be realized by the standard reduction procedure in this subsection.

\medskip
Multiplying (\ref{nonlinearprojection-helix}) by $ \frac{1}{(1+|V_d|^2)^2}\overline{Z_d} $ and integrating, we obtain
\begin{align*}
c\, \mbox{Re} \Bigg(\,\int_{\R^2} \frac{1}{(1+|V_d|^2)^2}\,Z_d\, {\overline {Z_d}}\,\Bigg)
\,=\,
&\,-\,\mbox{Re} \Bigg(\,\int_{\R^2} \frac{1}{(1+|V_d|^2)^2}\,  \overline{Z_d}\, {\mathcal E}\,\Bigg)
\\
&\,+\,\mbox{Re} \Bigg(\int_{\R^2} \frac{1}{(1+|V_d|^2)^2}\, \overline{Z_d}\,
\Big({\mathcal L}[\psi]+{\mathcal N}[\psi]+{\mathcal M}[\psi]\Big)\Bigg)
\\
\,=\,
&\,-\,\mbox{Re} \Bigg(\,\int_{\R^2} \frac{1}{(1+|V_d|^2)^2}\,  \overline{Z_d}\, {\mathbb E}_d\,\Bigg)
\\
&\,+\,\mbox{Re} \Bigg(\int_{\R^2} \frac{1}{(1+|V_d|^2)^2}\, \overline{Z_d}\,
\Big( {\mathbb L}_1[\phi]+{\mathbb N}_1[\phi]+{\mathbb M}_1[\phi]\Big)\Bigg).
\end{align*}
Using Proposition \ref{Proposition5.1} and the expression in (\ref{Nphi-helices}), we deduce that
\begin{equation*}
 \mbox{Re}  \Big(\int_{\R^2} \overline{Z_d}\, {\mathbb N}_1 [\phi] \Big) = o(\ve).
\end{equation*}
On the other hand, integration by parts, we have
 \[
 \mbox{Re}  \Big(\int_{\R^2} \frac{1}{(1+|V_d|^2)^2}\,  \overline{Z_d}\, {\mathbb L}_1 [\phi] \Big)
  = \mbox{Re}  \Big(\int_{\R^2} \frac{1}{(1+|V_d|^2)^2} \, \overline{\phi}\, {\mathbb L}_1 [Z_d] \Big).\]
Let us observe that
\begin{equation*}
\frac{\partial}{\partial d} {\mathbb S}_0 [V_d]
\,=\, {\mathbb L}_1 \Big[\frac{\partial V_d}{\partial d}\Big]
\,=\, {\mathbb L}_1 [Z_d] = O(\ve)
\end{equation*}
and thus by Proposition \ref{Proposition5.1}
\begin{equation*}
\mbox{Re}  \Big(\int_{\R^2} \frac{1}{(1+|V_d|^2)^2}\,  \overline{\phi}\, {\mathbb L}_1 [Z_d] \Big)
= o(\ve).
\end{equation*}

\medskip
As the strategy in standard reduction method, we are left to estimate the following integral
\begin{align}\label{errorprojection}
\,-\,\mbox{Re} \Bigg(\,\int_{\R^2} \frac{1}{(1+|V_d|^2)^2}\,  \overline{Z_d}\, {\mathbb E}_d\,\Bigg)
\,=\,
\mbox{Re} \Big(\int_{{\mathfrak S}_+ } \frac{1}{(1+|V_d|^2)^2}
\,\overline{Z_d}\, {\mathbb S}_4 [V_d]\, \Big)
\end{align}
The expressions of these error terms are given in (\ref{S4(Vd)}).
We will estimate the above integrals in the sequel.

\medskip
On ${\mathfrak S}_+, $ recall that $z= \vec{e}_1 +y$.
Then we have the estimates for the first term in (\ref{errorprojection})
\begin{align*}
\begin{aligned}
&\mathrm{Re}\Big(\int_{{\mathfrak S}_+ } \frac{1}{(1+|V_d|^2)^2}  \,\overline{Z_d}\,  {\mathbb S}_0 [V_d]\, \Big)
\\
& \,=\, -\,\int_{ {\mathfrak S}_+} \frac{1}{(1+ \rho^2)^2} \Bigg(\frac{1}{ d+ y_1} \frac{\partial \tilde{\rho}}{\partial y_1}  + \frac{\tilde{\rho} (\tilde{\rho}^2-1)}{\tilde{\rho}^2+1} \big(2 \nabla \varphi_s \nabla \theta +|\nabla \varphi_s|^2\big) \Bigg) \frac{\partial \tilde{\rho}}{\partial y_1}
\\
&\quad\,-\,\int_{ {\mathfrak S}_+} \frac{1}{(1+ \rho^2)^2} \frac{2 \tilde{\rho} (1-\tilde{\rho}^2)}{1+ \tilde{\rho}^2} \nabla \tilde{\rho} \cdot \nabla \varphi_s \frac{\partial \theta}{\partial y_1}
\,+\,O(\ve).
\end{aligned}
\end{align*}
Let us notice that
\[  \nabla \varphi_s= - \frac{1}{2 d}  \log d  \nabla y_2 + O(\ve \log |y| ).
\]
Hence
\begin{eqnarray*}
 \int_{{\mathfrak S}_+} \frac{1}{ (1+\rho^2)^2} \frac{1}{d + y_1} (\frac{\partial \tilde{\rho}}{\partial y_1} )^2\,=\, O(\ve),
\end{eqnarray*}
and also
\begin{eqnarray*}
& &  - \int_{ {\mathfrak S}_+} \frac{1}{(1+\rho^2)^2} \Biggl[ \frac{\,\tilde{\rho} (\tilde{\rho}^2-1)\,}{\tilde{\rho}^2+1} \big(2 \nabla \varphi_s \cdot \nabla \theta +|\nabla \varphi_s|^2\big) \frac{\partial \tilde{\rho}}{\partial y_1}
 + \frac{\,2 \tilde{\rho} (1-\tilde{\rho}^2)\,}{1+ \tilde{\rho}^2} \nabla \tilde{\rho} \cdot \nabla \varphi_s \frac{\partial \theta}{\partial y_1} \Biggl] +O(\ve)
\\
& & =  -\frac{1}{\,d \,} \log d  \int_{ {\mathfrak S}_+} \frac{  1-\rho^2 }{(\rho^2+1)^3 }  \frac{\rho \rho' }{|y|} dy \,+\, O(\ve)
\\
& & = -\frac{\pi }{\,8 d \,} \log d  \,+\, O(\ve).
\end{eqnarray*}

Recall the estimates in (\ref{IntegralQ2Pairs}) and (\ref{IntegralT2Pairs}).
the last terms  can be estimated by
\begin{align*}
\begin{aligned}
\kappa\, \varepsilon|\log\varepsilon|\,\mathrm{Re}\Big[\int_{{\mathfrak S}_+ } \frac{1}{(1+|V_d|^2)^2}
\,\overline{Z_d}\,{\mathbb T}_2 [V_d]\, \Big]
\,=\,-\kappa\, \varepsilon|\log\varepsilon|\frac{\pi}{2},
\\
\\
 \varepsilon|\log\varepsilon|\,\mathrm{Re}
\Big[\int_{{\mathfrak S}_+ } \frac{1}{(1+|V_d|^2)^2}
\,\overline{Z_d}\,{\mathbb Q}_2 [V_d]\, \Big]
\,=\,
 \varepsilon|\log\varepsilon|\frac{\pi}{4}.
\end{aligned}
\end{align*}

\medskip
As a conclusion, there holds
\begin{align}\label{reducealgebrahelix}
c(d)
\,=\,
c_0\,\pi\,\Bigg[\, -\frac{1 }{\,8 d\,} \log d
+\frac{\,1-2\kappa\,}{4}\ve |\log \ve | + O(\ve)\,\Bigg],
\end{align}
where $c_0 \not =0$.
Therefore, for any $2\kappa-1< 0$, we obtain a solution to $c(d)=0$ with the following asymptotic behavior:
\begin{align}
\frac{1 }{\,d\,} \log d
\sim 2(1-2\kappa)\ve |\log \ve|.
\end{align}

\begin{remark}\label{remarkLocationandveloHelix}
Recall the parameters $\ve$, $\kappa$,
given in (\ref{veHelix}), (\ref{kveHelix}).
By the relation in (\ref{reducealgebrahelix}),
we have that, for problem (\ref{Ishimoriproblem}), the vortex ring travels along $s_3$ axis
when its speed is sufficiently small with relation between it geometric parameters and traveling velocity
$$
\frac{1}{d}\log d\sim \frac{4c}{\sqrt{1-c^2}}-\frac{4\omega}{\sqrt{1-c^2}}.
$$
Similarly, for problem (\ref{Wavemaps}), there holds
$$
\frac{1}{d}\log d\sim \frac{4c}{\sqrt{1-c^2}}.
$$
\end{remark}

\bigskip
\noindent {\bf Acknowledgment.}
J. Yang is supported by the foundations:  NSFC(No.10901108), NSF of Guangdong (No.10451806001004770)
and the Foundation for Distinguished Young Talents in Higher Education of Guangdong (LYM11115).
Part of this work was done when J. Yang visited Chern Institute of Mathematics, Nankai University:
he is very grateful to the institution for the kind hospitality.

\end{document}